\documentclass [12pt]{amsart}
\usepackage{amssymb,amsmath,amsthm}
\newtheorem{theorem}{Theorem}
\newtheorem{lemma}{Lemma}
\newtheorem{proposition}{Proposition}

\newcommand{\ad}{\,\mathrm{ad}\,}

\newcommand{\GL}{\,\mathrm{GL}\,}
\newcommand{\SL}{\,\mathrm{SL}\,}

\newcommand{\diag}{\,\mathrm{diag}\,}
\begin{document}

\begin{center}

{\Large {\bf Automorphisms of Chevalley groups\\

\bigskip

of types $A_l, D_l, E_l$ over local rings without $1/2$\footnote{The
work is supported by the grant of Russian Fond of Basic Research 05-01-01048.} }}

\bigskip
\bigskip

{\large \bf E.~I.~Bunina}

\end{center}
\bigskip

\begin{center}

{\bf Abstract.}

\end{center}

In the given paper we prove that every automorphism of a Chevalley group of type $A_l$, $D_l$, or $E_l$, $l\geqslant 3$, over a commutative local ring without~$1/2$ is standard, i.\,e., it is a composition of ring, inner, central and graph  automorphisms.

\bigskip

\section*{Introduction}\leavevmode

Let $G_{\pi}$ be a Chevalley--Demazure group scheme associated with an irreducible root system~$\Phi$ of types $A_l$ ($l\geqslant 3$),
$D_l$ ($l\geqslant 4$), $E_l$ ($l=6,7,8$); $G_{\pi}(\Phi,R)$ be a set of points~$G_{\pi}$ with values
in~$R$; $E_{\pi}(\Phi,R)$ be the elementary subgroup of~$G_{\pi}(\Phi,R)$, where $R$
    is a commutative ring. In the given paper we describe automorphisms of the groups $G_{\pi}(\Phi,R)$ and $E_{\pi}(\Phi,R)$ over local commutative rings where $2$ is not invertible.

Similar results for Chevalley groups over fields were proved by
R.\,Steinberg~\cite{Stb1} for the finite case and by J.\,Humphreys~\cite{H} for the infinite one. Many papers were devoted
to description of automorphisms of Chevalley groups over different
commutative rings, we can mention here the papers of
Borel--Tits~\cite{v22}, Carter--Chen~Yu~\cite{v24},
Chen~Yu~\cite{v25}--\cite{v29}, E.\,Abe~\cite{Abe_OSN}, A.\,Klyachko~\cite{Klyachko}.

In the paper \cite{ravnyekorni}, it is shown that every automorphism of an elementary Chevalley group of type $A_l,D_l$, or $ E_l$, $l\geqslant 2$, over a local ring with $1/2$ is the composition of a ring automorphism and an
 \emph{automorphism--conjugation}, where an automorphism--conjugation is a conjugation of the Chevalley group in its adjoint representation by some matrix from the normalizer of this group in $\GL(V)$. In the paper \cite{normalizers} with the help of results of \cite{ravnyekorni} we proved that every automorphism of an arbitrary (elementary) Chevalley group of the type under consideration is standard, i.\,e., it is the composition of ring, inner, central and graph automorphisms. In the same paper we described the normalizers of Chevalley groups in the adjoint representation, this theorem holds also for the local rings without~$1/2$.

  In the papers \cite{bunF4}, \cite{korni2}, \cite{BunBl} by the similar methods we proved that automorphisms of Chevalley groups of types $F_4$, $G_2$, $B_l, l\geqslant 2$, over local rings with $1/2$ (for the type $G_2$ also with $1/3$)
are standard.

In the given paper we consider automorphisms of Chevalley groups of types $A_l,D_l, E_l$, $l\geqslant 3$, over local rings without $1/2$, and prove that every such automorphism is the composition of ring, inner, central and graph automorphisms. We use the results of the paper~\cite{normalizers} in our proof. The results of the given work can help to close the question on automorphisms of Chevalley groups over commutative rings without~$1/2$.

Note that while we consider root systems
$A_l$, $D_l$, $E_l$ in this paper, the case $A_l$ was completely closed by the papers of W.C.~Waterhouse~\cite{v46}, V.M.~Petechuk~\cite{v12},  Fuan Li and
Zunxian Li~\cite{v37}, and also for rings without~$1/2$.

The author is thankful to N.A.\,Vavilov,  A.A.\,Klyachko,
A.V.\,Mikhalev for valuable advices, remarks and discussions.

\section{Definitions and main theorems.}\leavevmode

We fix the root system~$\Phi$ of rank $>1$. Detailed texts about root
systems and their properties can be found in the books
\cite{Hamfris}, \cite{Burbaki}).
Suppose now that we have a
semisimple complex Lie algebra~$\mathcal L$ of type $\Phi$ with
Cartan subalgebra~$\mathcal H$ (detailed information about
semisimple Lie algebras can be found in the book~\cite{Hamfris}).

Then we can choose a basis $\{ h_1, \dots, h_l\}$ in~$\mathcal H$ and for every
$\alpha\in \Phi$ elements $x_\alpha \in {\mathcal L}_\alpha$ so that $\{ h_i; x_\alpha\}$ form a basis in~$\mathcal L$ and for every two elements of this basis their commutator is an integral linear combination of the elements of the same basis.

Let us introduce elementary Chevalley groups (see,
for example,~\cite{Steinberg}).

Let  $\mathcal L$ be a semisimple Lie algebra (over~$\mathbb C$)
with a root system~$\Phi$, $\pi: {\mathcal L}\to gl(V)$ be its
finitely dimensional faithful representation  (of dimension~$n$). If
$\mathcal H$ is a Cartan subalgebra of~$\mathcal L$, then a
functional
 $\lambda \in {\mathcal H}^*$ is called a
 \emph{weight} of  a given representation, if there exists a nonzero vector $v\in V$
 (that is called a  \emph{weight vector}) such that
for any $h\in {\mathcal H}$ $\pi(h) v=\lambda (h)v.$

In the space~$V$ there exists a basis of weight vectors such that
all operators $\pi(x_\alpha)^k/k!$ for $k\in \mathbb N$ are written
as integral (nilpotent) matrices. This basis is called a
\emph{Chevalley basis}. An integral matrix also can be considered as
a matrix over an arbitrary commutative ring with~$1$. Let $R$ be
such a ring. Consider matrices $n\times n$ over~$R$, matrices
$\pi(x_\alpha)^k/k!$ for
 $\alpha\in \Phi$, $k\in \mathbb N$ are included in $M_n(R)$.

Now consider automorphisms of the free module $R^n$ of the form
$$
\exp (tx_\alpha)=x_\alpha(t)=1+tx_\alpha+t^2 (x_\alpha)^2/2+\dots+
t^k (x_\alpha)^k/k!+\dots
$$
Since all matrices $x_\alpha$ are nilpotent, we have that this
series is finite. Automorphisms $x_\alpha(t)$ are called
\emph{elementary root elements}. The subgroup in $Aut(R^n)$,
generated by all $x_\alpha(t)$, $\alpha\in \Phi$, $t\in R$, is
called an \emph{elementary adjoint Chevalley group} (notation:
$E_{\ad}(\Phi,R)=E_{\ad}(R)$).

The action of $x_\alpha(t)$ on the Chevalley basis is described in
\cite{v23}, \cite{VavPlotk1}.

All weights of a given representation (by addition) generate a
lattice (free Abelian group, where every  $\mathbb Z$-basis  is also
a $\mathbb C$-basis in~${\mathcal H}^*$), that is called the
\emph{weight lattice} $\Lambda_\pi$.

Elementary Chevalley groups are defined not even by a representation
of the Chevalley groups, but just by its \emph{weight lattice}.
Namely, up to an abstract isomorphism an elementary Chevalley group
is completely defined by a root system~$\Phi$, a commutative
ring~$R$ with~$1$ and a weight lattice~$\Lambda_\pi$.

Among all lattices we can mark  the lattice corresponding to the
adjoint representation: it is generated by all roots (the \emph{root
lattice}~$\Lambda_{ad}$). The corresponding (elementary) Chevalley group is called \emph{adjoint}.

Introduce now Chevalley groups (see~\cite{Steinberg},
\cite{Chevalley}, \cite{v3}, \cite{v23}, \cite{v30}, \cite{v43},
\cite{VavPlotk1}, and also latter references in these papers).

Consider semisimple linear algebraic groups over algebraically
closed fields. These are precisely elementary Chevalley groups
$E_\pi(\Phi,K)$ (see.~\cite{Steinberg},~\S\,5).

All these groups are defined in $SL_n(K)$ as  common set of zeros of
polynomials of matrix entries $a_{ij}$ with integer coefficients
 (for example,
in the case of the root system $C_l$ and the universal
representation we have $n=2l$ and the polynomials from the condition
$(a_{ij})Q(a_{ji})-Q=0$). It is clear now that multiplication and
taking inverse element are also defined by polynomials with integer
coefficients. Therefore, these polynomials can be considered as
polynomials over arbitrary commutative ring with a unit. Let some
elementary Chevalley group $E$ over~$\mathbb C$ be defined in
$SL_n(\mathbb C)$ by polynomials $p_1(a_{ij}),\dots, p_m(a_{ij})$.
For a commutative ring~$R$ with a unit let us consider the group
$$
G(R)=\{ (a_{ij})\in \SL_n(R)\mid \widetilde p_1(a_{ij})=0,\dots
,\widetilde p_m(a_{ij})=0\},
$$
where  $\widetilde p_1(\dots),\dots \widetilde p_m(\dots)$ are
polynomials having the same coefficients as
$p_1(\dots),\dots,p_m(\dots)$, but considered over~$R$.

This group is called the \emph{Chevalley group} $G_\pi(\Phi,R)$ of
the type~$\Phi$ over the ring~$R$, and for every algebraically
closed field~$K$ it coincides with the elementary Chevalley group.

The subgroup of diagonal (in the standard basis of weight vectors)
matrices of the Chevalley group $G_\pi(\Phi,R)$ is called the
 \emph{standard maximal torus}
of $G_\pi(\Phi,R)$ and it is denoted by $T_\pi(\Phi,R)$. This group
is isomorphic to $Hom(\Lambda_\pi, R^*)$.

Let us denote by $h(\chi)$ the elements of the torus $T_\pi
(\Phi,R)$, corresponding to the homomorphism $\chi\in Hom
(\Lambda(\pi),R^*)$.

In particular, $h_\alpha(u)=h(\chi_{\alpha,u})$ ($u\in R^*$, $\alpha
\in \Phi$), where
$$
\chi_{\alpha,u}: \lambda\mapsto u^{\langle
\lambda,\alpha\rangle}\quad (\lambda\in \Lambda_\pi).
$$

Note that the condition
$$
G_\pi (\Phi,R)=E_\pi (\Phi,R)
$$
is not true even for fields, that are not algebraically closed.
Let us show the difference between Chevalley groups and their elementary subgroups in the case when $R$ is semilocal. In this case
 $G_\pi (\Phi,R)=E_\pi(\Phi,R)T_\pi(\Phi,R)$
(see~\cite{Abe1}), and elements $h(\chi)$ are connected with elementary generators
by the formula
\begin{equation}\label{ee4}
h(\chi)x_\beta (\xi)h(\chi)^{-1}=x_\beta (\chi(\beta)\xi).
\end{equation}

 In the case of semilocal rings from the formula~\eqref{ee4} we can see that
$$
[G(\Phi,R),G(\Phi,R)]\subseteq E(\Phi,R).
$$
In the case of the root systems $A_l,D_l, E_l$, $l\geqslant 2$, that are under consideration in this work, we have
$$x_{\alpha+\beta}(t)=[x_\alpha(t),x_\beta(1)], \quad \alpha+\beta\in \Phi,
$$
therefore
$$
[G(\Phi,R),G(\Phi,R)]=[E(\Phi,R), E(\Phi,R)]=E(\Phi,R).
$$

Define four types of automorphisms of a Chevalley group
 $G_\pi(\Phi,R)$, we
call them  \emph{standard}.

{\bf Central automorphisms.} Let $C_G(R)$ be a center of
$G_\pi(\Phi,R)$, $\tau: G_\pi(\Phi,R) \to C_G(R)$ be some
homomorphism of groups. Then the mapping $x\mapsto \tau(x)x$ from
$G_\pi(\Phi,R)$ onto itself is an automorphism of $G_\pi(\Phi,R)$,
that is denoted by~$\tau$ and called a \emph{central automorphism}
of the group~$G_\pi(\Phi,R)$.

{\bf Ring automorphisms.} Let $\rho: R\to R$ be an automorphism of
the ring~$R$. The mapping $(a_{i,j})\mapsto (\rho (a_{i,j}))$ from $G_\pi(\Phi,R)$
onto itself is an automorphism of the group $G_\pi(\Phi,R)$, that is
denoted by the same letter~$\rho$ and is called a \emph{ring
automorphism} of the group~$G_\pi(\Phi,R)$. Note that for all
$\alpha\in \Phi$ and $t\in R$ an element $x_\alpha(t)$ is mapped to
$x_\alpha(\rho(t))$.

{\bf Inner automorphisms.} Let $S$ be some ring containing~$R$,  $g$
be an element of $G_\pi(\Phi,S)$, that normalizes the subgroup $G_\pi(\Phi,R)$. Then
the mapping $x\mapsto gxg^{-1}$  is an automorphism
of the group~$G_\pi(\Phi,R)$, that is denoted by $i_g$ and is called an
\emph{inner automorphism}, \emph{induced by the element}~$g\in G_\pi(\Phi,S)$. If $g\in G_\pi(\Phi,R)$, then call $i_g$ a \emph{strictly inner}
automorphism.

{\bf Graph automorphisms.} Let $\delta$ be an automorphism of the
root system~$\Phi$ such that $\delta \Delta=\Delta$. Then there
exists a unique automorphisms of $G_\pi (\Phi,R)$ (we denote it by
the same letter~$\delta$) such that for every $\alpha \in \Phi$ and
$t\in R$ an element $x_\alpha (t)$ is mapped to
$x_{\delta(\alpha)}(\varepsilon(\alpha)t)$, where
$\varepsilon(\alpha)=\pm 1$ for all $\alpha \in \Phi$ and
$\varepsilon(\alpha)=1$ for all $\alpha\in \Delta$.

 Similarly we can define four type of automorphisms of the elementary
subgroup~$E(R)$. An automorphism~$\sigma$ of the group
 $G_\pi(\Phi,R)$ (or $E_\pi(\Phi,R)$)
is called  \emph{standard} if it is a composition of automorphisms
of these introduced four types.

Together with standard automorphisms we use the following ''temporary'' type of automorphisms of an elementary adjoint Chevalley group:

{\bf Automorphisms--conjugations.} Let $V$ be the representation space of the group $E_{\ad} (\Phi,R)$, $C\in \GL(V)$ be some matrix from its normalizer:
$$
C E_{\ad}(\Phi,R) C^{-1}= E_{\ad} (\Phi,R).
$$
 Then the mapping
 $x\mapsto CxC^{-1}$ from $E_\pi(\Phi,R)$ onto itself is an automorphism of the Chevalley group, that is denoted by~$i_С$ and is called an  \emph{automorphism--conjugation} of~$E(R)$,
\emph{induced by an element}~$C$ of~$\GL(V)$.

Our aim is to prove the next main theorem:

\begin{theorem}\label{main}
Let $G=G_{\pi}(\Phi,R)$ $(E_\pi(\Phi,R))$
be an \emph{(}elementary\emph{)} Chevalley group with a the root system $A_l, D_l$, or $E_l$, $l\geqslant 3$, $R$ is a commutative local ring where $2$ is not invertible. Then every automorphism of~$G$ is standard. If the Chevalley group is adjoint, then an inner automorphism in the composition is strictly
inner.
\end{theorem}

To prove this theorem we will start with the following important technical theorem:

\begin{theorem}\label{new}
Every automorphism of an elementary adjoint Chevalley group of a type
$A_l,D_l$, or $E_l$, $l\geqslant 2$, over a local ring without $1/2$ is a composition of a ring automorphism and
an automorphism--conjugation.
\end{theorem}

To apply Theorem~\ref{new}, we will use the following theorem from~\cite{normalizers}:

\begin{theorem}\label{norm}
Every automorphism--conjugation of an elementary adjoint Chevalley group of a type $A_l, D_l$, or $E_l$, $l\geqslant 3$, over a local
ring with $1$ is the composition of a strictly inner automorphism
\emph{(}conjugation with some element of the corresponding Chevalley group\emph{)} and a graph automorphism.
\end{theorem}

Three next sections are devoted to the proof of Theorem~\ref{new}.

An example of a Chevalley group of the type $A_2$ over a local ring without $1/2$, which has nonstandard automorphisms can be found
in~\cite{Petechuk2}.

\section{Replacing the initial automorphism to the special isomorphism.}\leavevmode

From this section we suppose that  $R$ is a local ring without $1/2$, the Chevalley group is adjoint, the root system is one of
the systems under consideration. In this section we use some reasonings
from~\cite{Petechuk1}.

Let $J$ be the maximal ideal (radical) of~$R$, $k$  the residue
field $R/J$. Then $E_J=E_{ad}( \Phi,R,J)$ is the greatest normal
proper subgroup of $E_{\ad}(\Phi,R)$ (see~\cite{Abe1}). Therefore,
$E_J$ is invariant under the action of~$\varphi$.

By this reason  the automorphism
$$
\varphi: E_{\ad} (\Phi,R)\to E_{\ad}(\Phi,R)
$$
induces an automorphism
$$
\overline \varphi: E_{\ad} (\Phi,R)/E_J=E_{\ad} (\Phi,k)\to
E_{\ad}(\Phi,k).
$$
The group $E_{\ad}(\Phi,k)$ is a Chevalley group over field,
therefore the automorphism $\overline \varphi$ is standard, i.\,e., it has the form
$$
\overline \varphi =    i_{\overline g} \overline \rho,\quad
\overline g\in N(E_{\ad}(\Phi,k)),
$$
where $\overline \rho$ is a ring automorphism, induced by some automorphism of~$k$.

It is clear that there exists a matrix $g\in GL_n(R)$ such that
its image under factorization  $R$ by~$J$ coincides with~$\overline
g$. We are not sure that  $g\in N_{\GL_n(R)}(E_{\ad}(\Phi,R))$.

Consider  a mapping $\varphi'= i_{g^{-1}} \varphi$. It is an
isomorphism of the group
 $E_{ad}(\Phi,R)\subset GL_n(R)$ onto some subgroup in $GL_n(R)$,
with the property that its image under factorization $R$ by $J$
coincides with the automorphism $\overline \rho$.

The above arguments prove

\begin{proposition}\label{dop1}
Every matrix $A\in E_{\ad}(\Phi,R)$ with elements from the subring~$R'$
of~$R$, generated by unit, is mapped under~$\varphi'$ to some matrix from
$$
A\cdot \GL_n(R,J)=\{ B\in \GL_n(R)\mid A-B\in M_n(J)\}.
$$
\end{proposition}

Let $a\in E_{\ad} (\Phi,R)$, $a^3=1$. Then the element $e=\frac{1}{3}
(1+a+a^2)$ is an idempotent from the ring $M_n(R)$. This idempotent
$e$ defines a decomposition of the free $R$-module $V\cong R^n$:
$$
V=eV\oplus (1-e)V=V_0\oplus V_1
$$
(the modules $V_0$, $V_1$ are free, since every projective module over a local ring is free~\cite{Mc}). Let $\overline V=\overline
V_0 \oplus \overline V_1$ be a decomposition of the $k$-module (linear space) $\overline V\cong k^n$ with respect to~$\overline a$, and
 $\overline e=\frac{1}{3} (1+\overline a+ \overline a^2)$.

Тогда имеем

\begin{proposition}\label{pr1_1}
The modules \emph{(}subspaces\emph{)}
 $\overline V_0$, $\overline V_1$ are images of the modules $V_0$, $V_1$ under factorization by~$J$.
\end{proposition}

\begin{proof} Let us denote the images of $V_0$, $V_1$ under factorization
by $J$ by $\widetilde V_0$, $\widetilde V_1$, respectively. Since
$V_0=\{ x\in V| ex=x\},$ $V_1= \{ x\in V|ex=0\},$  we have
 $\overline e(\overline x)=\frac{1}{2}(1+\overline a)(\overline x)=\frac{1}{2}
(1+\overline a(\overline
x))=\frac{1}{2}(1+\overline{a(x)})=\overline{e(x)}$. Then
$\widetilde V_0\subseteq \overline V_0$, $\widetilde V_1\subseteq
\overline V_1$.

Let $x=x_0+x_1$, $x_0\in V_0$, $x_1\in V_1$. Then $\overline
e(\overline x)=\overline e(\overline x_0)+\overline e (\overline
x_1)=\overline x_0$. If $\overline x\in \widetilde V_0$, then
$\overline x=\overline x_0$.
\end{proof}

{\bf Remark.}

Suppose now that $a^3=1$.
Let us study this matrix on the submodules  $V_0$ and $V_1$.

 Let$x\in V_1$. Then $x=ex=\frac{1}{3}(a^2+a+1)(x)$, therefore $(a^2+a-2)(x)=0$. So we have $(a-1)(a^2+a-2)(x)=0$, since $a^3=1$, then $(a-1)(a^2+a+1)(x)=0$. Subtracting one equation from another one we get $3(a-1)(x)=0$ and as $3$ is invertible, we have $a(x)=x$. Consequently, on the submodule~$V_1$ the mapping $a$ is identical.

Now let $x\in V_0$. Let $a(x)=\lambda x$. Then it is clear that $\lambda^2+\lambda+1=0$. Since $\lambda \ne 1$, then $\lambda$ is a root from the unit of the third power, not equal to~$1$ (there are two of them, since $3$ is invertible). If $y=a(x)\ne \lambda x$, then the vectors $x$ and $y$ are linearly independent, also $a(y)=a^2(x)=-x-a(x)=-x-y$. Therefore, on the submodule  $\langle x,y\rangle$ the matrix $a$ acts invariantly and has the form
\begin{equation}\label{bl}
\begin{pmatrix}
0& -1\\
1& -1
\end{pmatrix}.
\end{equation}
Thus the module $V_0$ is a direct sum of invariant submodules of dimensions one or two, and on $1$-dimensional submodules the matrix $a$ acts on the vectors $x$ as the multiplication to $\xi$ or $\xi^2$ ($\xi^3=1$).

Clear that if the initial matrix $a$ has integer coefficients, then its matrix invariants (trace, determinant, etc.) are also intefer numbers, so they are integer in any basis. Suppose that the number of $1$-dimensional submodules for $\xi$ is~$p$, the number of them for $\xi^2$ is~$q$. Suppose that $p\geqslant q$. Then determinant on the part of basis generated by these $1$-dimensional submodules is $\xi^{p-q}$. If this number is integer, then  $3|(p-q)$. The trace of the matrix on the given basis part is $-q+(p-q)\xi$, i.\,e., $p-q$ is zero in the ring~$R$. Since all another matrix invariant are also integer, and also all odd numbers in~$R$ are invertible, we have that $p=q$.
Therefore if a matrix $a$ has integer coefficients, then in some basis it consists of a unit block of the dimension $\dim V_1$ and of $2\times 2$ blocks of the form \eqref{bl} and of complete dimension $\dim V_0$.

\medskip

 Let now for a matrix~$a$ with integer coefficients
$b=\varphi'(a)$. Then $b^3=1$ and $b$ is equivalent to $a$ modulo~$J$.

\begin{proposition}\label{pr1_2}
Suppose that $ a,b\in E_\pi(\Phi,R)$, $a^3=b^3=1$, $a$ is a matrix with elements from the subring $R'\subset R$, generated by unit, $b$ and $a$ are equivalent modulo~$J$, $V=V_0\oplus V_1$ is the decomposition of~$V$ with respect to~$a$, $V=V_0'\oplus V_1'$ is the decomposition of~$V$ with respect to~$b$. Then $\dim V_0'=\dim V_0$, $\dim
V_1'=\dim V_1$.
\end{proposition}

\begin{proof}
We have an  $R$-basis of the module~$V$ $\{ e_1,\dots,e_n\}$ such that $\{
e_1,\dots,e_k\}\subset V_0$, $\{ e_{k+1},\dots, e_n\}\subset V_1$.
Clear that
$$
\overline a \overline e_i=\overline{ae_i}=\overline {(\sum_{j=1}^n
a_{ij} e_j)}= \sum_{j=1}^n \overline a_{ij} \overline e_j.
$$
Let $\overline V=\overline V_0\oplus \overline V_1$, $\overline V
= \overline V_0'\oplus \overline V_1'$ be decompositions of $k$-module (linear space)~$\overline V$ with respect to $\overline a$ and $
\overline b$. Clear that $\overline V_0= \overline V_0'$, $\overline
V_1 =\overline V_1'$. Therefore, by Proposition~\ref{pr1_1}
the images of the modules $V_0$ and $V_0'$, $V_1$ and $V_1'$ under factorization by~$J$
coincide. Let us take such $\{ f_1,\dots, f_k\}\subset V_0'$, $\{
f_{k+1},\dots, f_n\}\subset V_1'$, that $\overline f_i=\overline
e_i$, $i=1,\dots,n$. Since the basis change matrix (from basis $\{ e_1,\dots,
e_n\}$ to $\{ f_1,\dots, f_n\}$) is invertible (is equivalent to the unit matrix modulo~$J$), then $\{ f_1,\dots, f_n\}$ is an $R$-basis
in~$V$. Clear that $\{ f_1,\dots, f_k\}$ is an $R$-basis in
$V_0'$, $\{ v_{k+1},\dots, v_n\}$ is an $R$-basis in~$V_1'$.
\end{proof}

From this proposition, the above remark and equivalence of the matrices  $a$ and $b$ it follows that for $b$ there exists some basis of the module~$V$,
in which $b$ has the same form, as  $a$ in the initial basis. Therefore, $a$ and $b$ are conjugate.

\section{The images of $w_{\alpha_i}x_{\alpha_i}(1)$ and some elements of the Weil group}

Consider some fixed elementary adjoint Chevalley group
$E=E_{\ad}(\Phi,R)$ with a root system $A_l$ ($l\geqslant 3$), $D_l$ ($l\geqslant
4$), $E_6$, $E_7$, or $E_8$, its adjoint representation in the group $\GL_n(R)$ ($n=l+2m$, where $m$ is the number of positive roots of~$\Phi$), with weight vector basis $v_1=x_{\alpha_1},
v_{-1}=x_{-\alpha_1}, \dots, v_n=x_{\alpha_n}, v_{-n}=x_{-\alpha_n},
V_1=h_{1},\dots,V_l=h_{l}$, corresponding to the Chevalley basis of~$\Phi$.

We also have an isomorphism~$\varphi'$, described in Section~2.

In the beginning we consider the matrix $Q_1=w_{\alpha_1}x_{\alpha_1}(1)$ in our basis. Note that $Q_1^3=E$.

On the basis part generated by $\{ v_1, v_{-1}, V_{1}, V_2\}$, this matrix is invariant and has the form
$$
\begin{pmatrix}
0& -1&  0& 0\\
-1& 1&  2& -1\\
0& -1&  -1& 1\\
0& 0&  0& 1
\end{pmatrix};
$$
 on the basis part generated by all $v_j, v_{-j}$, where $\langle
\alpha_1, \alpha_j\rangle=0$, and also by all corresponding $V_i$, $i>2$,
it is identical; on the basis part $\{ v_j,v_{-j}, v_k, v_{-k}\}$,
where $\langle \alpha_1,\alpha_j\rangle =-1$,
$\alpha_k=\alpha_1+\alpha_j$, it has the form
$$
\begin{pmatrix}
-1& 0& 1& 0\\
0& 0& 0& 1\\
-1& 0& 0& 0\\
0& -1& 0& -1 \end{pmatrix}.
$$

Note that the given matrix has only integer coefficients, therefore its image under the isomorphism is conjugate to it.

Now add  (temporarily) to the ring $R$
an element~$\xi$ such that $\xi^3=1$. In other words we take such an outer element and consider the ring $\overline R$, generated by~$R$ and~$\xi$.

Then in some basis the matrix  $Q_1$ has a diagonal form with $1,\xi,\xi^2$ on the diagonal. Namely,
the basis part of the first form is mapped to $\diag[\xi,\xi^2,1,1]$, the basis part of the second type is mapped to  $\diag[1,1]$, the basis part of the third type is mapped to $\diag[\xi,\xi,\xi^2,\xi^2]$.
Clear that other elements $Q_i=w_{\alpha_i}x_{\alpha_i}(1)$ have similar properties.

For different root systems from the list under consideration we take now the following sets of roots:

--- for the root system $A_l$ we consider the set $\alpha_1=e_1-e_2$, $\alpha_3=e_3-e_4$, \dots, $\alpha_{l-1}=e_{l-1}-e_l$, or $\alpha_l=e_l-e_{l+1}$, depending on evenness of~$l$, i.\,e., just take simple roots next nearest;

--- for the root system  $D_l$ we consider the set $\alpha_1=e_1-e_2$, $\alpha_3=e_3-e_4$, \dots, $\alpha_{l-1}=e_{l-1}-e_l$, $e_1+e_2$, $e_3+e_4$, \dots, $e_{l-1}+e_l$;

--- for the system $E_8$ ($E_6$, $E_7$ are similar) consider the set $\alpha_1=e_1-e_2$, $\alpha_4=e_3-e_4$, $\alpha_6=e_5-e_6$, $\alpha_8=e_7-e_8$, $\frac{1}{2}(e_1+e_2-e_3-e_4-e_5-e_6-e_7-e_8)$, $\frac{1}{2}(e_1+e_2+e_3+e_4-e_5-e_6-e_7-e_8)$, $\frac{1}{2}(e_1+e_2+e_3+e_4+e_5+e_6-e_7-e_8)$, $\frac{1}{2}(e_1+e_2+e_3+e_4+e_5+e_6+e_7+e_8)$.

Let us denote the obtained set (sequence) of rots by $\gamma_1,\dots, \gamma_k$.

Note that all roots $\gamma_1,\dots,\gamma_k$ are mutually orthogonal. It means that all matrices $Q_{\gamma_1},\dots, Q_{\gamma_k}$ mutually commute, and so their images $P_{\gamma_i}=\varphi'(Q_{\gamma_i})$ commute. Therefore, we can in the same basis write them in diagonal form with the same (the same number and the same places) $1,\xi,\xi^2$ on the diagonal. Temporarily we come to these basis.

Consider the transition matrix of this basis.

On the basis part $\{ \gamma_i, -\gamma_i, h_{\gamma_i}, h_{\alpha}\}$, where $\langle \alpha,\gamma_i\rangle =-1$,
the transition matrix has the form
$$
\begin{pmatrix}
1& -\xi& -2\xi^2& -\xi\\
-\xi& 1& 2\xi^2& \xi\\
1& -1& 1& 0\\
0& 0& 0& 1
\end{pmatrix}.
$$

If we take such a root $\alpha$, that $\langle \alpha,\gamma_i\rangle=\langle \alpha,\gamma_j\rangle=-1$, then on the basis part, generated by the roots
$$
\{ \alpha, -\alpha, \alpha+\gamma_i, -\alpha-\gamma_i, \alpha+\gamma_j, -\alpha-\gamma_j,
\alpha+\gamma_i+\gamma_j, -\alpha-\gamma_i-\gamma_j\},
$$
the transition matrix is
$$
\begin{pmatrix}
1& 0& \xi& 0& \xi& 0& -\xi& 0\\
0& 1& 0& -\xi& 0& -\xi& 0& -\xi\\
\xi& 0& 1& 0& -\xi& 0& \xi& 0\\
0& -\xi& 0& 1& 0& -\xi& 0& -\xi\\
\xi& 0& -\xi& 0& 1& 0& \xi& 0\\
0& -\xi& 0& -\xi& 0& 1& 0& -\xi\\
-\xi& 0& \xi& 0& \xi& 0& 1& 0\\
0& -\xi& 0& -\xi& 0& -\xi& 0& 1
\end{pmatrix}.
$$

Note also that for any pair of roots $\gamma_i,\gamma_j$, $1\leqslant i,j \leqslant k$, there exists a root $\gamma_{i,j}$ such that $\langle \gamma_i, \gamma_{i,j}\rangle=\langle \gamma_j, \gamma_{i,j}\rangle =-1$ (for the root pair $e_p-e_{p+1}$ and $e_q\pm e_{q+1}$ it is the root $e_{p+1}-e_q$; for the root pair $e_p+e_{p+1}$ and $e_q\pm e_{q+1}$ it is the root $-e_{p+1}-e_q$; for the root pair $e_p-e_{p+1}$ and $e_p+e_{p+1}$ it is the root $-e_p+e_q$, $q\ne p,p+1$; for the root pair $e_p-e_{p+1}$ and $\frac{1}{2}(e_1+e_2\pm\dots \pm (e_p+e_{e_{p+1}})\pm \dots)$ it is the root $e_1\mp e_p$, or $e_1\mp e_{p+1}$; finally, for the pair $\frac{1}{2}(e_1+e_2-e_3-e_4-e_5-e_6-e_7-e_8)$ and $\frac{1}{2}(e_1+e_2+e_3+e_4-e_5-e_6-e_7-e_8)$ it is $e_1-e_5$).

Consider the Weil group element $w_{i,j}=w_{\gamma_{i,j}}(1)w_{\gamma_i}(1) w_{\gamma_j}(1) w_{\gamma_{i,j}}(1)$. It is easy to show that $w_{i,j}^2=E$ and $w_{i,j}Q_{\gamma_i}w_{i,j}=Q_{\gamma_j}$.

The basis change under consideration does not change the elements $w_{i,j}$.

Look to $W_{i,j}=\varphi_1(w_{i,j})$. Clear that it is also an element of order two, that maps  $P_{\gamma_i}$ to $P_{\gamma_j}$ by conjugation (now we are in the basis, where all $P_q$ are diagonal and coincide with the diagonal form of~$Q_q$).
Also $W_{i,j}$ commutes with all $P_{\gamma_k}$, $k\ne i,j$.

Consider some roots $\alpha,\beta$ and the place $(\alpha,\beta)$ in the matrix $W_{i,j}$.

Let the roots $\alpha$ and $\beta$ be orthogonal to $\gamma_i$ and $\gamma_j$. Then they are necessary not orthogonal to some  $\gamma_k$, $k\ne i,j$. Clear that to have not zero on the place $(\alpha,\beta)$, it is necessary that for all $\gamma_k$, $k\ne i,j$, $\langle \alpha, \gamma_k\rangle=\langle \beta, \gamma_k\rangle$.

For the root system $A_l$ we can suppose for our convenience that $\gamma_i=\alpha_1$, $\gamma_j=\alpha_3$, $\alpha=e_p-e_q$, $\beta=e_t-e_s$. So we directly obtain $p,q,t,s>4$. Let $p\ne t$, $q\ne s$. Consider $\gamma_k=e_p-e_{p+1}$ or $e_{p-1}-e_p$. Clear that by our supposition in the first case we have $s=p+1$, in the second case we have $s=p-1$. Similarly, $t=q+1$ or $q-1$. We cannot have the situation $p\ne t$, $q=s$, or $p=t$, $q\ne s$. Therefore to have something not equal to zero on the place $(\alpha,\beta)$ in the case under consideration, it is necessary $\alpha=\beta$, or $\alpha=-\beta\pm \gamma_l\pm \gamma_k$ (for every $\alpha$ there exist not more than one such $\beta$).

For the root system $D_l, E_l$ the situation is just better, since in the sequence $\gamma_1,\dots,\gamma_k$ there are additional roots to differ $\alpha$ and $\beta$. Consequently in this case we have $\alpha=\beta$.

Now let the root  $\alpha$ be orthogonal to $\gamma_i$ and $\gamma_j$, and $\beta$ be not orthogonal (and we suppose that it is not orthogonal at least to $\gamma_j$).
We know that  $P_{\gamma_i}W_{i,j}=W_{i,j}P_{\gamma_j}$. Let on the place $(\alpha,\beta)$ in the matrix $W_{i,j}$ be~$a$, then in the condition in the left matrix on this place is still~$a$, and in the right one it is  $a\xi$, or $a\xi^2$.   Therefore $a=0$.

Now we only need to consider the cases when both roots  $\alpha,\beta$ are not orthogonal to some of roots $\gamma_i,\gamma_j$. From the same arguments as above it follows that if both $\alpha$ and $\beta$ are orthogonal to one root of $\gamma_i,\gamma_j$ and are not orthogonal to the second one, then on the place $(\alpha,\beta)$ it is zero.

Suppose now, for definiteness, that $\alpha \perp \gamma_j$, $\beta\perp \gamma_i$. In this case we again have two variants: either $\alpha=\pm \gamma_i$, or there exist also some $\gamma_t$, $t\ne i,j$, which are not orthogonal to $\alpha$. In the first case it is clear that the coefficient is not zero only if $\beta=\pm \gamma_j$ (with the same sign). In the second case we again consider the root system $A_l$, for the convenience we suppose that $\gamma_i=\alpha_1$, $\gamma_j=\alpha_3$, then $\alpha=\pm e_1\pm e_p$ or $\pm e_2\pm e_p$, $p>4$, and we see that for $\beta$ there are two possibilities: either it is the root $w_{i,j}(\alpha)=\pm e_3\pm e_p$ or $\pm e_4\pm e_p$, or it is the root $-w_{i,j}(\alpha)\pm \gamma_j\pm \gamma_t$, $\gamma_t=e_p-e_{p+1}$ or $e_{p-1}-e_p$. For other root systems the situation is again just more defined.

In the last case $\alpha$ and $\beta$ are not orthogonal neither to $\gamma_i$, no to $\gamma_j$.
Let us again consider root systems separately.

For the root system $A_l$ in supposition $\gamma_i=\alpha_1$, $\gamma_j=\alpha_3$, we have $\alpha=\pm e_p \pm e_q$, where $p=1,2$, $q=3,4$. Without loss of generality we can suppose that $\alpha=e_1-e_3$, then $\beta=e_3-e_1$, or $\beta=e_2-e_4$. Therefore  for every possible root $\alpha$ there are two possibilities for~$\beta$.

For the root system $D_l$ if $\gamma_i=e_1-e_2$, $\gamma_j=e_3-e_4$, $\alpha=e_1-e_3$, $\beta=e_3-e_1$, then $\langle \alpha,e_1+e_2\rangle \ne \langle \beta, e_1+e_2\rangle$, therefore we again have zero on the place $(\alpha,\beta)$. Consequently, we have now only possibility  $\beta=e_2-e_4$. Similarly we obtain the unique possibility  $\beta=w_{i,j}(\alpha)$ for other choices of $\alpha$. If $\gamma_i=e_1-e_2$, $\gamma_j=e_1+e_2$, then $\alpha$ is necessary  $\pm e_1\pm e_p$, or $\pm e_2\pm e_p$, $p>2$. Let for example $\alpha=e_1-e_3$. Then we directly see $\beta=e_1-e_3$, only one possibility.

For the root system $E_l$ the situation is completely the same as for~$D_l$.

On the basis part generated by $h_{\alpha_1},\dots, h_{\alpha_l}$, the matrix $W_{i,j}$ can be arbitrary yet.

Note that for the root system $D_l,E_l$ the elements $W_{i,j}$ are directly (simply after taking the conditions) close to their inverse images $w_{i,j}$, therefore we will consider basis changes only for the root systems $A_l$ as for the case, where $W_{i,j}$ are the most far from $w_{i,j}$. For convenience we consider the root system~$A_5$. Clear that for all other root systems under consideration changes will be similar.

In the root system $A_5$ we order the roots as follows: $\pm(e_1-e_2), \pm (e_2-e_3), \pm (e_3-e_4), \pm (e_4-e_5), \pm (e_5-e_6)$,
$\pm(e_1-e_3), \pm (e_2-e_4), \pm (e_3-e_5), \pm (e_4-e_6)$, $\pm (e_1-e_4), \pm(e_2-e_5), \pm (e_2-e_6)$, $\pm (e_1-e_5), \pm (e_2-e_6),
\pm (e_1-e_6)$.

We are interested in $W_{1,3}$ and $W_{3,5}$ (since $W_{1,5}$ is generated by them). Let us consider separately basis parts, invariant under both $W_{1,3}$ and $W_{3,5}$.

The following basis part is invariant under the both matrices
$
\{ \pm(e_1-e_2), \pm (e_3-e_4), \pm (e_5-e_6)\},
$
and
{\footnotesize
$$
W_{1,3}=
\begin{pmatrix}
0& 0& a_{1,5}& 0& 0& 0\\
0& 0& 0& a_{2,6}& 0& 0\\
a_{5,1}& 0& 0& 0& 0& 0\\
0& a_{6,2}& 0& 0& 0& 0\\
0& 0& 0& 0& a_{9,9}& 0\\
0& 0& 0& 0& 0& a_{10,10}
\end{pmatrix},\quad
W_{3,5}=
\begin{pmatrix}
b_{1,1}& 0& 0& 0& 0& 0\\
0& b_{2,2}& 0& 0& 0& 0\\
0& 0& 0& 0& b_{5,9}& 0\\
0& 0& 0& 0& 0& b_{6,10}\\
0& 0& b_{9,5}&  0& 0& 0\\
0& 0& 0& b_{10,6}& 0& 0
\end{pmatrix}.
$$
}

Let us make the basis change (commuting with all $P_{\alpha_1}, P_{\alpha_3}, P_{\alpha_5}$) as follows:
$v_1'=v_1, v_2'=v_2, v_5'=a_{1,5}v_5, v_6'=a_{2,6}v_6, v_9'=b_{5,9}a_{1,5}v_9, v_{10}'=b_{6,10}a_{1,6}v_{10}$. Then the matrices $W_{1,3}, W_{3,5}$ on these basis parts are
$$
\begin{pmatrix}
0& 0& 1& 0& 0& 0\\
0& 0& 0& 1& 0& 0\\
a_{5,1}'& 0& 0& 0& 0& 0\\
0& a_{6,2}'& 0& 0& 0& 0\\
0& 0& 0& 0& a_{9,9}& 0\\
0& 0& 0& 0& 0& a_{10,10}
\end{pmatrix}
\text{ и }
\begin{pmatrix}
b_{1,1}& 0& 0& 0& 0& 0\\
0& b_{2,2}& 0& 0& 0& 0\\
0& 0& 0& 0& 1& 0\\
0& 0& 0& 0& 0& 1\\
0& 0& b_{9,5}'&  0& 0& 0\\
0& 0& 0& b_{10,6}'& 0& 0
\end{pmatrix},
$$
and since $W_{1,3}^2=W_{3,5}^2=E$, then $a_{5,1}'=a_{6,2}'=b_{9,5}'=b_{10,6}'=1$. Becides,
$$
W_{1,3}W_{3,5}=\begin{pmatrix}
0& 0& 0& 0& 1& 0\\
0& 0& 0& 0& 0& 1\\
b_{1,1}& 0& 0& 0& 0& 0\\
0& b_{2,2}& 0& 0& 0& 0\\
0& 0& a_{9,9}& 0& 0& 0\\
0& 0& 0& a_{10,10}& 0& 0
\end{pmatrix}
$$
is an element of the order three, so $b_{1,1}=a_{9,9}$, $b_{2,2}=a_{10,10}$, all these elements has the order~$2$.

The next basis part is
$$
\pm (e_2-e_3), \pm(e_4-e_5), \pm (e_1-e_4), \pm (e_2-e_5), \pm(e_3-e_6),  \pm (e_1-e_6).
$$

The matrix $W_{1,3}$ on it is equal to
{\footnotesize
$$
\left(\begin{array}{cccccccccccc}
a_{3,3}& 0& 0& 0& 0& a_{3,20}& 0& 0& 0& 0& 0& 0\\
0& a_{4,4}& 0& 0& a_{4,19}& 0& 0& 0& 0& 0& 0& 0\\
0& 0& 0& 0& 0& 0& a_{7,21}& 0& 0& 0& 0& a_{7,30}\\
0& 0& 0& 0& 0& 0& 0& a_{8,22}& 0& 0& a_{8,29}& 0\\
0& a_{19,4}& 0& 0& a_{19,19}& 0& 0& 0& 0& 0& 0& 0\\
a_{20,3}& 0& 0& 0& 0& a_{20,20}& 0& 0& 0& 0& 0& 0\\
0& 0& a_{21,7}& 0& 0& 0& 0& 0& 0& a_{21,24}& 0& 0\\
0& 0& 0& a_{22,8}& 0& 0& 0& 0& a_{22,23}& 0& 0& 0\\
0& 0& 0& 0& 0& 0& 0& a_{23,22}& 0& 0& a_{23,29}& 0\\
0& 0& 0& 0& 0& 0& a_{24,21}& 0& 0& 0& 0& a_{24,30}\\
0& 0& 0& a_{29,8}& 0& 0& 0& 0& a_{29,23}& 0& 0& 0\\
0& 0& a_{30,7}& 0& 0& 0& 0& 0& 0& a_{30,24}& 0& 0
\end{array}\right),
$$
}
and the matrix $W_{3,5}$ is
{\footnotesize
$$
\left(\begin{array}{cccccccccccc}
0& 0& 0& 0& 0& 0& b_{3,21}& 0& 0& 0& 0& b_{3,30}\\
0& 0& 0& 0& 0& 0& 0& b_{4,22}& 0& 0& b_{4,29}& 0\\
0& 0& b_{7,7}& 0& 0& 0& 0& 0& 0& b_{7,24}& 0& 0\\
0& 0& 0& b_{8,8}& 0& 0& 0& 0& b_{8,23}& 0& 0& 0\\
0& 0& 0& 0& 0& 0& 0& b_{19,22}& 0& 0& b_{19,29}& 0\\
0& 0& 0& 0& 0& 0& b_{20,21}& 0& 0& 0& 0& b_{20,30}\\
b_{21,3}& 0& 0& 0& 0& b_{21,20}& 0& 0& 0& 0& 0& 0\\
0& b_{22,4}& 0& 0& b_{22,19}& 0& 0& 0& 0& 0& 0& 0\\
0& 0& 0& b_{23,8}& 0& 0& 0& 0& b_{23,23}& 0& 0& 0\\
0& 0& b_{24,7}& 0& 0& 0& 0& 0& 0& b_{24,24}& 0& 0\\
0& b_{29,4}& 0& 0& b_{29,19}& 0& 0& 0& 0& 0& 0& 0\\
b_{30,3}& 0& 0& 0& 0& b_{30,20}& 0& 0& 0& 0& 0& 0
\end{array}\right).
$$
}

Let us describe on this basis part the required basis change (certainly, it commutes with all $P_{\alpha_1}, P_{\alpha_3}, P_{\alpha_5}$). We will make this change by steps.

1. All basis vectors are not moved, except the vectors $v_{19}$ and $v_{20}$. Namely, $v_{19}'=a_{4,4}v_4+a_{19,4}v_{19}$, $v_{20}'=a_{3,3}v_3+a_{20,3}v_{20}$. After such transformations the structure of $W_{3,5}$ is not changed (only nonzero elements $b_{i,j}$ can change values, but for simplicity we will not write any primes), and the matrix $W_{1,3}$ (also since it has the order two) will have the form
{
$$
\left(\begin{array}{cccccccccccc}
0& 0& 0& 0& 0& 1& 0& 0& 0& 0& 0& 0\\
0& 0& 0& 0& 1& 0& 0& 0& 0& 0& 0& 0\\
0& 0& 0& 0& 0& 0& a_{7,21}& 0& 0& 0& 0& a_{7,30}\\
0& 0& 0& 0& 0& 0& 0& a_{8,22}& 0& 0& a_{8,29}& 0\\
0& 1& 0& 0& 0& 0& 0& 0& 0& 0& 0& 0\\
1& 0& 0& 0& 0& 0& 0& 0& 0& 0& 0& 0\\
0& 0& a_{21,7}& 0& 0& 0& 0& 0& 0& a_{21,24}& 0& 0\\
0& 0& 0& a_{22,8}& 0& 0& 0& 0& a_{22,23}& 0& 0& 0\\
0& 0& 0& 0& 0& 0& 0& a_{23,22}& 0& 0& a_{23,29}& 0\\
0& 0& 0& 0& 0& 0& a_{24,21}& 0& 0& 0& 0& a_{24,30}\\
0& 0& 0& a_{29,8}& 0& 0& 0& 0& a_{29,23}& 0& 0& 0\\
0& 0& a_{30,7}& 0& 0& 0& 0& 0& 0& a_{30,24}& 0& 0
\end{array}\right).
$$
}
Now let us make the following change.

2. All basis vectors are remained, except  $v_{21}$ and $v_{22}$. Namely, $v_{21}'=b_{21,3}v_{21}+b_{30,3}v_{30}$, $v_{22}'=b_{22,4}v_{22}+b_{29,4}v_{29}$. After such transformations the structure of $W_{1,3}$ is not changed, and $W_{3,5}$ is
$$
\left(\begin{array}{cccccccccccc}
0& 0& 0& 0& 0& 0& 1& 0& 0& 0& 0& 0\\
0& 0& 0& 0& 0& 0& 0& 1& 0& 0& 0& 0\\
0& 0& b_{7,7}& 0& 0& 0& 0& 0& 0& b_{7,24}& 0& 0\\
0& 0& 0& b_{8,8}& 0& 0& 0& 0& b_{8,23}& 0& 0& 0\\
0& 0& 0& 0& 0& 0& 0& b_{19,22}& 0& 0& b_{19,29}& 0\\
0& 0& 0& 0& 0& 0& b_{20,21}& 0& 0& 0& 0& b_{20,30}\\
1& 0& 0& 0& 0& 0& 0& 0& 0& 0& 0& 0\\
0& 1& 0& 0& 0& 0& 0& 0& 0& 0& 0& 0\\
0& 0& 0& b_{23,8}& 0& 0& 0& 0& b_{23,23}& 0& 0& 0\\
0& 0& b_{24,7}& 0& 0& 0& 0& 0& 0& b_{24,24}& 0& 0\\
0& 0& 0& 0& b_{29,19}& 0& 0& 0& 0& 0& 0& 0\\
0& 0& 0& 0& 0& b_{30,20}& 0& 0& 0& 0& 0& 0
\end{array}\right).
$$

3. All basis vectors are not changed, except the vectors $v_{7}$ and $v_{8}$. Namely, $v_{7}'=a_{7,21}v_{7}+a_{24,21}v_{24}$, $v_{8}'=a_{8,22}v_{8}+a_{23,22}v_{23}$. After such transformations the structure of $W_{3,5}$ is not changed, and the matrix $W_{1,3}$ is
$$
\left(\begin{array}{cccccccccccc}
0& 0& 0& 0& 0& 1& 0& 0& 0& 0& 0& 0\\
0& 0& 0& 0& 1& 0& 0& 0& 0& 0& 0& 0\\
0& 0& 0& 0& 0& 0& 1& 0& 0& 0& 0& a_{7,30}\\
0& 0& 0& 0& 0& 0& 0& 1& 0& 0& a_{8,29}& 0\\
0& 1& 0& 0& 0& 0& 0& 0& 0& 0& 0& 0\\
1& 0& 0& 0& 0& 0& 0& 0& 0& 0& 0& 0\\
0& 0& 1& 0& 0& 0& 0& 0& 0& a_{21,24}& 0& 0\\
0& 0& 0& 1& 0& 0& 0& 0& a_{22,23}& 0& 0& 0\\
0& 0& 0& 0& 0& 0& 0& 0& 0& 0& a_{23,29}& 0\\
0& 0& 0& 0& 0& 0& 0& 0& 0& 0& 0& a_{24,30}\\
0& 0& 0& 0& 0& 0& 0& 0& a_{29,23}& 0& 0& 0\\
0& 0& 0& 0& 0& 0& 0& 0& 0& a_{30,24}& 0& 0
\end{array}\right).
$$

 We know that $W_{3,5}^2=E$. The matrix $W_{3,5}^2$ on the place $(5,2)$ has $b_{19,22}$, on the place  $(6,1)$ it has $b_{20,21}$. Therefore, $b_{19,22}=b_{20,21}=0$.

Now we will make the later changes.

4. All basis vectors are not changed, except the vectors $v_{23}$ and $v_{24}$. Namely, $v_{23}'=b_{7,7}v_{7}+b_{24,7}v_{24}$, $v_{24}'=b_{8,8}v_{8}+b_{23,8}v_{24}$. After these transformations the structure if $W_{1,3}$ is not changed, and the matrix $W_{3,5}$ is now
$$
\left(\begin{array}{cccccccccccc}
0& 0& 0& 0& 0& 0& 1& 0& 0& 0& 0& 0\\
0& 0& 0& 0& 0& 0& 0& 1& 0& 0& 0& 0\\
0& 0& 0& 0& 0& 0& 0& 0& 0& 1& 0& 0\\
0& 0& 0& 0& 0& 0& 0& 0& 1& 0& 0& 0\\
0& 0& 0& 0& 0& 0& 0& 0& 0& 0& b_{19,29}& 0\\
0& 0& 0& 0& 0& 0& 0& 0& 0& 0& 0& b_{20,30}\\
1& 0& 0& 0& 0& 0& 0& 0& 0& 0& 0& 0\\
0& 1& 0& 0& 0& 0& 0& 0& 0& 0& 0& 0\\
0& 0& 0& 1& 0& 0& 0& 0& 0& 0& 0& 0\\
0& 0& 1& 0& 0& 0& 0& 0& 0& 0& 0& 0\\
0& 0& 0& 0& b_{29,19}& 0& 0& 0& 0& 0& 0& 0\\
0& 0& 0& 0& 0& b_{30,20}& 0& 0& 0& 0& 0& 0
\end{array}\right).
$$

Since the matrix $W_{1,3}W_{3,5}$ has the order~$3$, we directly obtain $a_{22,23}=a_{21,24}=0$. Since $W_{1,3}$ has the order~$2$, we have $a_{7,30}=a_{8,29}=0$.

5. Now let us consider the last basis change, where $v_{29}'=b_{29,19}v_{29}$, $v_{30}'=b_{30,20}v_{30}$, $v_{23}'=a_{23,29}v_{23}$, $v_{24}'=a_{24,30}v_{30}$.
After it we have $W_{1,3}=w_{1,3}$, $W_{3,5}=w_{3,5}$ on the basis part under consideration.

On the basis part $\pm(e_1-e_3), \pm (e_2-e_4), \pm (e_3-e_5), \pm (e_4-e_6), \pm (e_1-e_5), \pm (e_2-e_6)$ all arguments are similar to the previous part.

We now need only to consider the basis part $h_{\alpha_1}, h_{\alpha_2}, h_{\alpha_3}, h_{\alpha_4}, h_{\alpha_5}$. We can make an arbitrary basis change there, since all $P_{\alpha_1}, P_{\alpha_3}, P_{\alpha_5}$ on this part are identical.

Suppose that on this part
$$
W_{1,3}=\begin{pmatrix}
a_{31,31}& a_{31,32}& a_{31,33}& a_{31,34}& a_{31,35}\\
a_{32,31}& a_{32,32}& a_{32,33}& a_{32,34}& a_{32,35}\\
a_{33,31}& a_{33,32}& a_{33,33}& a_{33,34}& a_{33,35}\\
a_{34,31}& a_{34,32}& a_{34,33}& a_{34,34}& a_{34,35}\\
a_{35,31}& a_{35,32}& a_{35,33}& a_{35,34}& a_{35,35}
\end{pmatrix},\quad
W_{3,5}=\begin{pmatrix}
b_{31,31}& b_{31,32}& b_{31,33}& b_{31,34}& b_{31,35}\\
b_{32,31}& b_{32,32}& b_{32,33}& b_{32,34}& b_{32,35}\\
b_{33,31}& b_{33,32}& b_{33,33}& b_{33,34}& b_{33,35}\\
b_{34,31}& b_{34,32}& b_{34,33}& b_{34,34}& b_{34,35}\\
b_{35,31}& b_{35,32}& b_{35,33}& b_{35,34}& b_{35,35}
\end{pmatrix}.
$$

At first let us make the following basis change. All basis vectors are not moved, only $v_{33}'=a_{31,33}v_{31}+a_{32,33}v_{32}+a_{33,33}v_{33}+a_{34,33}v_{34}+a_{35,33}v_{35}$. After such a change the matrix $W_{1,3}$ has the form
$$
W_{1,3}=\begin{pmatrix}
0& a_{31,32}& 1& a_{31,34}& a_{31,35}\\
0& a_{32,32}& 0& a_{32,34}& a_{32,35}\\
1& a_{33,32}& 0& a_{33,34}& a_{33,35}\\
0& a_{34,32}& 0& a_{34,34}& a_{34,35}\\
0& a_{35,32}& 0& a_{35,34}& a_{35,35}
\end{pmatrix}.
$$

By the next basis change we use only $v_{32}$. Namely, $v_{32}'=(a_{31,34}v_{31}+a_{32,34}v_{32}+a_{33,34}v_{33}+a_{34,34}v_{34}+a_{35,34}v_{35})-v_{34}-v_{33}$. Then the matrix $W_{1,3}$ will have the form  (we use also the condition $W_{1,3}^2=E$)
$$
W_{1,3}=\begin{pmatrix}
0& -1& 1& 0& a_{31,35}\\
0& -1& 0& 1& a_{32,35}\\
1& -1& 0& 1& a_{33,35}\\
0& 0& 0& 1& a_{34,35}\\
0& 0& 0& 0& a_{35,35}
\end{pmatrix}.
$$

By similar changes for $W_{3,5}$ (we change the basis elements $v_{34}$ and $v_{35}$) we can obtain
$$
W_{3,5}=\begin{pmatrix}
b_{31,31}& 0& 0& 0& 0\\
b_{32,31}& 1& 0& 0& 0\\
b_{33,31}& 1& 0& -1& 0\\
b_{34,31}& 1& 0& -1& 0\\
b_{35,31}& 0& 1& -1& 1
\end{pmatrix}.
$$

Then we have variables  $a_{31,35}, a_{32,35}, a_{33,35}, a_{34,35}, a_{35,35}, b_{31,31}, b_{32,31}, b_{33,31}, b_{34,31}, b_{35,31}$, where $a_{35,35}$ and $b_{31,31}$ are equivalent to~$1$ modulo radical, all other variables are form the radical.

Besides, we have three matrix conditions:

1) $W_{1,3}^2=E$;

2) $W_{3,5}^2=E$;

3) $(W_{1,3}W_{3,5})^3=E$.

So we obtain some polynomial conditions on the variables.

We will remove variables by steps, considering their coefficients modulo radical. If finally we obtain that the last variable is zero, then all variables are equals to integer numbers, which are equivalent them modulo radical.

From the position $(1,5)$ of the condition 1) it follows $a_{32,35}=a_{33,35}+a_{31,35}$; from the position $(3,5)$ it follows $a_{34,35}=0$. Similarly from the position $(5,1)$ of the condition 2) it follows $b_{34,31}=b_{35,31}+b_{33,31}$; from the position $(33,31)$ it follows $b_{32,31}=0$. From the position $(2,2)$ of the condition 3) it follows $b_{35,31}=-b_{33,31}$; from the position $(1,3)$ it follows $a_{31,35}=0$; from the position $(1,4)$ it follows $b_{33,31}=0$; from the position $(1,2)$ it follows $a_{35,35}=1$; from the position $(1,1)$ it follows $b_{31,31}=1$, and, finally, from the position $(2,3)$ it follows $a_{33,35}=0$. Therefore, on this basis part also $W_{1,3}$ and $W_{3,5}$ coincide with $w_{1,3}$ and $w_{3,5}$, respectively.

In a whole we obtain that $W_{1,3}$ differs from $w_{1,3}$ just in the positions $(9,9)$ and $(10,10)$, where the matrix is diagonal, but on the diagonal it has not units but some elements of the order two, equivalent to the unit. Since $W_{1,3}$ has determinant~$1$ (because it is a product of commutators), these elements on the diagonal are equal. Similarly the elements on the positions $(1,1)$ and $(2,2)$ of the matrix $W_{3,5}$ are equal (an equal to multipliers for $W_{1,3}$). Let us denote them by $\mu$. Clear, how by diagonal basis change in all its parts under consideration we can obtain $W_{1,3}=\mu w_{1,3}$, $W_{3,5}=\mu w_{3,5}$.

Now by the basis change inverse to the initial one, we can return $Q_{\gamma_i}$ from the diagonal form to its normal form, where in the matrix there are no elements $\xi$ and $\xi^2$. As we know, $W_{i,j}$ are not changed.

Therefore we can now suppose that we have an isomorphism  $\varphi_2$ with all properties of  $\varphi_1$, and such that $\varphi_2(Q_{\gamma_i})=Q_{\gamma_i}$
for all $i=1,\dots,k$; $\varphi_2(w_{i,j})=\mu w_{i,j}$ for all $i,j=1,\dots,k$, $i\ne j$, $\mu^2=1$.

Below we suppose that we consider an isomorphism~$\varphi_2$ with this property.

\section{Limitation of the images of $x_{\alpha}(1)$ and $w_{\alpha}(1)$  on different basis parts.}

Suppose now that $\varphi_2(x_{\alpha_i}(1))=x_i$, $\varphi_2(w_{\alpha_i}(1))=W_i$.

At first we consider $x_1$ (we now suppose that all roots are numerated so that $\gamma_1=\alpha_1$). We know that $x_1$ commutes with all $Q_{\gamma_i}$, $i>1$. According to it $x_1$ is divided into some blocks. Let us see to these blocks.

Suppose now that we deal with the root system $A_l$, $l$ is odd.

Consider the root $\alpha=e_i-e_j$, $i<j$. If $\alpha=\alpha_1=e_1-e_2$, then $\alpha$ is orthogonal to all $\gamma_j$, $j>1$. Clear that if we take another arbitrary root~$\beta$, not collinear to~$\alpha$, then there exists some $\gamma_j$, not orthogonal to~$\beta$. Therefore on the place $(\alpha,\beta)$ in the matrix $x_1$ there is zero. So we see that the roots $\pm \alpha_1$ together with the basic elements $h_1,\dots,h_l$ gives us the separated invariant basis part.
Let now $\alpha=e_1-e_i$, $i>2$ (similarly we can consider the roots $e_2-e_i$, $i>2$). Such a root is not orthogonal to exactly one of the roots from the sequence $\gamma_2,\dots, \gamma_k$, for example, $\gamma_j$ (it is either $e_i-e_{i+1}$, or $e_{i-1}-e_i$).
The roots $\pm (e_1-e_i), \pm (e_2-e_i)$ have the same property and either $\pm (e_1-e_{i-1}), \pm (e_2-e_{i-1}), \pm (e_{i-1}-e_i)$, or $\pm (e_1-e_{i+1}), \pm (e_2-e_{i+1}), \pm (e_{i}-e_{i+1})$. Such set of roots gives a basis part, invariant under the matrix~$x_1$. Now we just need to consider the root $e_i-e_j$, which is nor orthogonal to the pair $\gamma_p,\gamma_q$. Without loss of generality suppose that $i,j$ are both even. Then the basis part $\pm (e_i-e_j), \pm(e_i-e_{j+1}), \pm (e_{i+1}-e_j), \pm (e_{i+1}-e_{j+1})$ is invariant.

Suppose now that the root system is again $A_l$, but $l$ is even.

By similar arguments we can see that the basis parts for the matrix $x_1$, have one of the following forms:

1) $\pm (e_1-e_2), \pm (e_1-e_{l+1}), \pm (e_2-e_{l+1}), h_1,\dots$;

2) $\pm (e_1-e_{2i-1}), \pm (e_1-e_{2i}), \pm (e_2-e_{2i-1}), \pm (e_2-e_{2i}), \pm (e_{2i-1}-e_{2i}), \pm (e_{2i-1}-e_{l+1}),
\pm (e_{2i}-e_{l+1})$, $i> 1$;

3) $\pm (e_{2i-1}-e_{2j-1}), \pm (e_{2i-1}-e_{2j}), \pm (e_{2i}-e_{2j-1}), \pm (e_{2i}-e_{2j})$, $i\ne j$, $i,j,>1$.

Now let us deal with the root system $D_l$.
The basis parts for this system and for the matrix $x_1$, are the following:

1) $\pm (e_1-e_2), h_1,\dots$;

2) $\pm (e_1+e_2)$;

3) $\pm (e_1-e_i), \pm (e_2-e_i)$, $i> 2$;

4) $\pm (e_1+e_i), \pm (e_2+e_i)$, $i> 2$;

5) $\pm (e_i-e_j)$, $i,j>2$;

6) $\pm (e_i+e_j)$, $i,j>2$.

We see that they are ``smaller'' than basis parts for~$A_l$.

Finally let us consider the root system $E_8$.
In this case we have the following basis parts:

1) $\pm (e_1-e_2), h_1,\dots$;

2) $\pm (e_1+e_2)$;

3)  $\pm (e_1-e_i), \pm (e_2-e_i)$, $i> 2$;

4)  $\pm (e_1+e_i), \pm (e_2+e_i)$, $i> 2$;

5) $\pm (e_i-e_j)$, $i,j>2$;

6) $\pm (e_i+e_j)$, $i,j>2$;

7) $\pm \frac{1}{2}(-e_1+e_2\pm e_3\pm e_4\pm e_5\pm e_6\pm e_7\pm e_8), \pm \frac{1}{2} (e_1-e_2\pm e_3\pm e_4\pm e_5\pm e_6\pm e_7\pm e_8)$;

8) $\pm \frac{1}{2}(e_1+e_2\pm e_3\pm e_4\pm e_5\pm e_6\pm e_7\pm e_8)$.

For such root system also the basis parts are strictly smaller than for the root systems~$A_l$.

According to it we can consider only the root systems $A_l$, $l\geqslant 3$.

Let us now look to the element~$x_2$. It commutes with $Q_i$ of our sequence, except the first two of them, therefore the invariant basis parts are bigger then for~$x_1$.

As we have seen above we do not need to consider root systems  $D_l$, $E_l$, we consider only the root systems $A_l$, $l\geqslant 3$.

For these systems if $l$ is odd, then $x_2$ is separated to the following parts:

1) $\pm (e_1-e_2), \pm (e_1-e_3), \pm (e_1-e_4), \pm (e_2-e_3), \pm (e_2-e_4), \pm (e_3-e_4), h_1,\dots$;

2) $\pm (e_1 - e_i), \pm (e_2-e_i), \pm (e_3-e_i), \pm (e_4-e_i), \pm (e_1-e_{i+1}), \pm (e_2-e_{i+1}) , \pm (e_3-e_{i+1}),
\pm (e_4-e_{i+1})$, $i>4$, $i$ is odd;

3) $\pm (e_i-e_j), \pm (e_i-e_{j+1}), \pm (e_{i+1}-e_j), \pm (e_{i+1}-e_{j+1})$, $i,j> 4$, $i,j$ are odd;

4) $\pm (e_i-e_{i+1})$, $i> 4$, $i$ is odd.

If  $l$ is even, then the basis parts are of the following form:

1) $\pm (e_1-e_2), \pm (e_1-e_3), \pm (e_1-e_4), \pm (e_2-e_3), \pm (e_2-e_4), \pm (e_3-e_4), \pm (e_1-e_{l+1}),
\pm (e_2-e_{l+1}), \pm (e_3-e_{l+1}), \pm (e_4-e_{l+1}), h_1,\dots$;

2) $\pm (e_1 - e_i), \pm (e_2-e_i), \pm (e_3-e_i), \pm (e_4-e_i), \pm (e_1-e_{i+1}), \pm (e_2-e_{i+1}) , \pm (e_3-e_{i+1}),
\pm (e_4-e_{i+1}), \pm (e_i-e_{l+1}), \pm (e_{i+1}-e_{l+1})$, $i>4$, $i$ is odd;

3) $\pm (e_i-e_j), \pm (e_i-e_{j+1}), \pm (e_{i+1}-e_j), \pm (e_{i+1}-e_{j+1})$, $i,j> 4$, $i,j$ are odd;

4) $\pm (e_i-e_{i+1}), \pm (e_i-e_{l+1}), \pm (e_{i+1}-e_{l+1})$, $i> 4$, $i$ is odd.

Therefore, in a whole  (for $x_1$ and $x_2$) our matrices are divided into the same basis parts, which a listed above for the matrix~$x_2$.

It is clear that the cases of even and odd~$l$ should be considered separately, but we do not need to consider any systems except~$A_l$.

\section{Images of $w_{\alpha_i}$ and $x_{\alpha_i}(1)$.}

Note that on every basis part situation is as follows: there are some known matrices (for example, $Q_{\alpha_1}$, $Q_{\alpha_3}$, $w_{1,3}$) and some unknown matrices (for example, the images of $w_{\alpha_1}(1)$ and $w_{\alpha_2}(1)$, they will be denoted by $w_{\alpha_1}(1)+W_1$ and $w_{\alpha_2}(1)+W_2$). Clear that both matrices $W_1$ and $W_2$ are from the ideal $M_N(J)$. All other unknown matrices are expressed by the known matrices and introduced above unknown matrices (for example, $\varphi_2(x_{\alpha_1}(1))=Q_{\alpha_1}\cdot (w_{\alpha_1}(1)+W_1)^3$). Besides, there is some set of conditions  (for example,
$(w_{\alpha_1}(1)+W_1)Q_{\alpha_3}=Q_{\alpha_3}(w_{\alpha_1}(1)+W_1)$), which are true for our unknown matrices $W_1$ and $W_2$. We want to show that it is possible after some additional basis changes (commuting with all known matrices) come to the situation when all the conditions are satisfied only for the zero matrices $W_1$ and $W_2$. It means that the matrices $w_{\alpha_1}(1)$ and $w_{\alpha_2}(1)$ under the obtained isomorphism are mapped to themselves, what is required. 

Let the elements of $W_1$ and $W_2$ be denoted by $z_1,\dots, z_p$.
Note that every matrix condition gives $N^2$ polynomial equations of variables $z_1,\dots,z_p$ with integer coefficients.

Suppose that one of these polynomials can be represented in the form
$$
z_{k_0}A+z_1B_1+\dots+z_{k_0-1}B_{k_0-1}+z_{k_0+1}B_{k_0+1}+\dots+z_{p}B_{p}=0,
$$
and the polynomial  $A$ is invertible modulo radical, $B_i$ are some polynomials ($z_{k_0}$ can enter in all polynomials, including~~$A$). Then
$$
z_{k_0}=-\frac{z_1B_1+\dots+z_{k_0-1}B_{k_0-1}+z_{k_0+1}B_{k_0+1}+\dots+z_{p}B_{p}}{A},
$$
we can substitute the expression for $z_{k_0}$ in all other polynomial conditions. If we can choose sequentially $p$ such condition, then in the process of the described substitution every time we remove one new variable, then the last condition will look like
$$
z_{k_{p}}C=0,
$$
where $C$ is some rational expression of variables
$z_1,\dots, z_{p}$, invertible modulo radical. Therefore we can say that $z_{k_{p}}=0$, and therefore all other variables are zeros. The existence of these $p$
conditions is equivalent to the existence of such $p$ conditions that the square matrix with entries equal to the coefficients of these conditions modulo radical, has an invertible (i.\,e., odd) determinant.

Since it is vary complicated to write this  $p\times p$ matrix, we can sequentially write the obtained condition, but for simplicity write the coefficients $A$ and $B_i$ modulo radical (in the result we will write just numbers $0$ and $1$).

Clear that this procedure is equivalent to the  ``linearization'' of all conditions up to the variables of the matrices $W_1$ and $W_2$. Namely, if after opening parentheses in the condition somewhere there is an expression $W_iW_jA$, where $i,j\in \{1,2\}$, $A$ is an arbitrary matrix, then such an expression we can suppose zero. Finally all condition becomes linear with respect to $W_1$ and $W_2$.

In the result we need to show that if we express one variables by another in these (linearized) conditions all variables become zeros.

At first we will show it on the simple basis parts: on the third and fourth one.

\subsection{Basis parts of the fourth type.}

It is the simplest basis of the form $\pm (e_i-e_{i+1})$, $i> 4$, $i$ is odd. On this basis part all matrices interesting for us ($w_{\alpha_1}(1)$, $w_{\alpha_2}(1)$, $x_{\alpha_1}(1)$, $x_{\alpha_2}(1)$, $x_{\alpha_1+\alpha_2}(1)$) are identical.

We know that the element $(w_{\alpha_1}(1)+W_1)(w_{\alpha_2}(1)+W_2)=(E+W_1)(E+W_2)$ has the order three. Linearizing this condition we obtain:
$$
(E+W_1+W_2)^3=E\Longleftrightarrow 3W_1=-3W_2 \Longleftrightarrow W_2=W_1.
$$

Therefore we can suppose $W_2=W_1$.

Also we have the condition $\varphi(x_{\alpha_1}(1))=x_{\alpha_1}(1)+X_1=E+X_1=Q_1\cdot (w_{\alpha_1}(1)+W_1)^3=(E+W_1)^3$, therefore after linearization $X_1=3W_1=W_1$. Besides, $E+X_{1+2}=\varphi_2(x_{\alpha_1+\alpha_2}(1))=(E+W_2)\cdot (E+X_1)\cdot (E+W_2)^3=E+X_1+4W_1=E+X_1$, i.\,e., $X_{1+2}=X_1$. Similarly, $E+X_2=\varphi_2(x_{\alpha_2}(1)=(E+W_1)\cdot (E+X_{1+2})\cdot (E+W_1)^3$, so$X_2=X_{1+2}=X_1=W_1$.

Now let us use the last condition $x_{\alpha_1+\alpha_2}(1)x_{\alpha_1}(1)x_{\alpha_2}(1)=x_{\alpha_2}(1)x_{\alpha_1}(1)$, which gives for the images
$$
(E+W_1)(E+W_1)(E+W_1)=(E+W_1)(E+W_1).
$$
Clear that it gives us $W_1=0$, what was required.

\subsection{Basis parts of the third form.}

Now we consider basis parts of the third type, namely,
$\pm (e_i-e_j), \pm (e_i-e_{j+1}), \pm (e_{i+1}-e_j), \pm (e_{i+1}-e_{j+1})$, $i,j> 4$, $i,j$ нечетны.

Note that on this part, as on the previous one, all matrices $w_{\alpha_1}(1)$, $w_{\alpha_2}(1)$, $x_{\alpha_1}(1)$, $x_{\alpha_2}(1)$, $x_{\alpha_1+\alpha_2}(1)$ are identical.

Therefore, all arguments are completely similar to the previous case, since we use there only the identical forms of matrices and two conditions, which hold only on the present basis part.

\subsection{Basis parts of the second form.}

Now we come to the basis part of the second part, namely,  $\pm (e_1 - e_i), \pm (e_2-e_i), \pm (e_3-e_i), \pm (e_4-e_i), \pm (e_1-e_{i+1}), \pm (e_2-e_{i+1}) , \pm (e_3-e_{i+1}),
\pm (e_4-e_{i+1})$, $i>4$, $i$ is odd. The matrix $W_1$ is decomposed into the direct sum on the basis parts $\pm (e_1 - e_i), \pm (e_2-e_i), \pm (e_1-e_{i+1}), \pm (e_2-e_{i+1})$ and $ \pm (e_3-e_i), \pm (e_4-e_i),  \pm (e_3-e_{i+1}),
\pm (e_4-e_{i+1})$, so the half of coefficients of~$W_1$ are exactly zeros.

The matrix $Q_{\alpha_i}$ is
{\scriptsize $$
Q_i=\left(\begin{array}{cccccccccccccccc}
-1& 0& 0& 0& 0& 0& 0& 0& 1& 0& 0& 0& 0& 0& 0& 0\\
0 & 0& 0& 0& 0& 0& 0& 0& 0& 1& 0& 0& 0& 0& 0& 0\\
0& 0& -1& 0& 0& 0& 0& 0& 0& 0& 1& 0& 0& 0& 0& 0\\
0& 0& 0& 0& 0& 0& 0& 0& 0& 0& 0& 1& 0& 0& 0& 0\\
0& 0& 0& 0& -1& 0& 0& 0& 0& 0& 0& 0& 1& 0& 0& 0\\
0& 0& 0& 0& 0& 0& 0& 0& 0& 0& 0& 0& 0& 1& 0& 0\\
0& 0& 0& 0& 0& 0& -1& 0& 0& 0& 0& 0& 0& 0& 1& 0\\
0& 0& 0& 0& 0& 0& 0& 0& 0& 0& 0& 0& 0& 0& 0& 1\\
-1& 0& 0& 0& 0& 0& 0& 0& 0& 0& 0& 0& 0& 0& 0& 0\\
0& -1& 0& 0& 0& 0& 0& 0& 0& -1& 0& 0& 0& 0& 0& 0\\
0& 0& -1& 0& 0& 0& 0& 0& 0& 0& 0& 0& 0& 0& 0& 0\\
0& 0& 0& -1& 0& 0& 0& 0& 0& 0& 0& -1& 0& 0& 0& 0\\
0& 0& 0& 0& -1& 0& 0& 0& 0& 0& 0& 0& 0& 0& 0& 0\\
0& 0& 0& 0& 0& -1& 0& 0& 0& 0& 0& 0& 0& -1& 0& 0\\
0& 0& 0& 0& 0& 0& -1& 0& 0& 0& 0& 0& 0& 0& 0& 0\\
0& 0& 0& 0& 0& 0& 0& -1& 0& 0& 0& 0& 0& 0& 0& -1
\end{array}\right);
$$}
the matrix $Q_{\alpha_1}$ is
{\scriptsize $$
Q_1=\left(\begin{array}{cccccccccccccccc}
-1& 0& 1& 0& 0& 0& 0& 0& 0& 0& 0& 0& 0& 0& 0& 0\\
0& 0& 0& 1& 0& 0& 0& 0& 0& 0& 0& 0& 0& 0& 0& 0\\
-1& 0& 0& 0& 0& 0& 0& 0& 0& 0& 0& 0& 0& 0& 0& 0\\
0& -1& 0& -1& 0& 0& 0& 0& 0& 0& 0& 0& 0& 0& 0& 0\\
0& 0& 0& 0& 1& 0& 0& 0& 0& 0& 0& 0& 0& 0& 0& 0\\
0& 0& 0& 0& 0& 1& 0& 0& 0& 0& 0& 0& 0& 0& 0& 0\\
0& 0& 0& 0& 0& 0& 1& 0& 0& 0& 0& 0& 0& 0& 0& 0\\
0& 0& 0& 0& 0& 0& 0& 1& 0& 0& 0& 0& 0& 0& 0& 0\\
0& 0& 0& 0& 0& 0& 0& 0& -1& 0& 1& 0& 0& 0& 0& 0\\
0& 0& 0& 0& 0& 0& 0& 0& 0& 0& 0& 1& 0& 0& 0& 0\\
0& 0& 0& 0 &0 &0 &0 &0& -1& 0& 0& 0& 0& 0& 0& 0\\
0& 0& 0& 0& 0& 0& 0& 0& 0& -1& 0& -1& 0& 0& 0& 0\\
0& 0& 0& 0& 0& 0& 0& 0& 0& 0& 0& 0& 1& 0& 0& 0\\
0& 0& 0& 0& 0& 0& 0& 0& 0& 0& 0& 0& 0& 1& 0& 0\\
0& 0& 0& 0& 0& 0& 0& 0& 0& 0& 0& 0& 0& 0& 1& 0\\
0& 0& 0& 0& 0& 0& 0& 0& 0& 0& 0& 0& 0& 0& 0& 1
\end{array}\right);
$$}
the matrix $w_{1,3}$ is
{\scriptsize $$
\left(\begin{array}{cccccccccccccccc}
0& 0& 0& 0& 1& 0& 0& 0& 0& 0& 0& 0& 0& 0& 0& 0\\
0& 0& 0& 0& 0& 1& 0& 0& 0& 0& 0& 0& 0& 0& 0& 0\\
0& 0& 0& 0& 0& 0& 1& 0& 0& 0& 0& 0& 0& 0& 0& 0\\
0& 0& 0& 0& 0& 0& 0& 1& 0& 0& 0& 0& 0& 0& 0& 0\\
1& 0& 0& 0& 0& 0& 0& 0& 0& 0& 0& 0& 0& 0& 0& 0\\
0& 1& 0& 0& 0& 0& 0& 0& 0& 0& 0& 0& 0& 0& 0& 0\\
0& 0& 1& 0& 0& 0& 0& 0& 0& 0& 0& 0& 0& 0& 0& 0\\
0& 0& 0& 1& 0& 0& 0& 0& 0& 0& 0& 0& 0& 0& 0& 0\\
0& 0& 0& 0& 0& 0& 0& 0& 0& 0& 0& 0& 1& 0& 0& 0\\
0& 0& 0& 0& 0& 0& 0& 0& 0& 0& 0& 0& 0& 1& 0& 0\\
0& 0& 0& 0& 0& 0& 0& 0& 0& 0& 0& 0& 0& 0& 1& 0\\
0& 0& 0& 0& 0& 0& 0& 0& 0& 0& 0& 0& 0& 0& 0& 1\\
0& 0& 0& 0& 0& 0& 0& 0& 1& 0& 0& 0& 0& 0& 0& 0\\
0& 0& 0& 0& 0& 0& 0& 0& 0& 1& 0& 0& 0& 0& 0& 0\\
0& 0& 0& 0& 0& 0& 0& 0& 0& 0& 1& 0& 0& 0& 0& 0\\
0& 0& 0& 0& 0& 0& 0& 0& 0& 0& 0& 1& 0& 0 & 0& 0
\end{array}\right);
$$}
the matrix $w_{\alpha_1}(1)$ is
{\scriptsize $$
w_1=\left(\begin{array}{cccccccccccccccc}
0& 0& 1& 0& 0& 0& 0& 0& 0& 0& 0& 0& 0& 0& 0& 0\\
0& 0& 0& 1& 0& 0& 0& 0& 0& 0& 0& 0& 0& 0& 0& 0\\
-1& 0& 0& 0& 0& 0& 0& 0& 0& 0& 0& 0& 0& 0& 0& 0\\
0& -1& 0& 0& 0& 0& 0& 0& 0& 0& 0& 0& 0& 0& 0& 0\\
0& 0& 0& 0& 1& 0& 0& 0& 0& 0& 0& 0& 0& 0& 0& 0\\
0& 0& 0& 0& 0& 1& 0& 0& 0& 0& 0& 0& 0& 0& 0& 0\\
0& 0& 0& 0& 0& 0& 1& 0& 0& 0& 0& 0& 0& 0& 0& 0\\
0& 0& 0& 0& 0& 0& 0& 1& 0& 0& 0& 0& 0& 0& 0& 0\\
0& 0& 0& 0& 0& 0& 0& 0& 0& 0& 1& 0& 0& 0& 0& 0\\
0& 0& 0& 0& 0& 0& 0& 0& 0& 0& 0& 1& 0& 0& 0& 0\\
0& 0& 0& 0& 0& 0& 0& 0& -1& 0& 0& 0& 0& 0& 0& 0\\
0& 0& 0& 0& 0& 0& 0& 0& 0& -1& 0& 0& 0& 0& 0& 0\\
0& 0& 0& 0& 0& 0& 0& 0& 0& 0& 0& 0& 1& 0& 0& 0\\
0& 0& 0& 0& 0& 0& 0& 0& 0& 0& 0& 0& 0& 1& 0& 0\\
0& 0& 0& 0& 0& 0& 0& 0& 0& 0& 0& 0& 0& 0& 1& 0\\
0& 0& 0& 0& 0& 0& 0& 0& 0& 0& 0& 0& 0& 0& 0& 1
\end{array}\right);
$$}
the matrix $w_{\alpha_2}(1)$ is
{\scriptsize $$
w_2=\left(\begin{array}{cccccccccccccccc}
1& 0& 0& 0& 0& 0& 0& 0& 0& 0& 0& 0& 0& 0& 0& 0\\
0& 1& 0& 0& 0& 0& 0& 0& 0& 0& 0& 0& 0& 0& 0& 0\\
0& 0& 0& 0& 1& 0& 0& 0& 0& 0& 0& 0& 0& 0& 0& 0\\
0& 0& 0& 0& 0& 1& 0& 0& 0& 0& 0& 0& 0& 0& 0& 0\\
0& 0& -1& 0& 0& 0& 0& 0& 0& 0& 0& 0& 0& 0& 0& 0\\
0& 0& 0& -1& 0& 0& 0& 0& 0& 0& 0& 0& 0& 0& 0& 0\\
0& 0& 0& 0& 0& 0& 1& 0& 0& 0& 0& 0& 0& 0& 0& 0\\
0& 0& 0& 0& 0& 0& 0& 1& 0& 0& 0& 0& 0& 0& 0& 0\\
0& 0& 0& 0& 0& 0& 0& 0& 1& 0& 0& 0& 0& 0& 0& 0\\
0& 0& 0& 0& 0& 0& 0& 0& 0& 1& 0& 0& 0& 0& 0& 0\\
0& 0& 0& 0& 0& 0& 0& 0& 0& 0& 0& 0& 1& 0& 0& 0\\
0& 0& 0& 0& 0& 0& 0& 0& 0& 0& 0& 0& 0& 1& 0& 0\\
0& 0& 0& 0& 0& 0& 0& 0& 0& 0& -1& 0& 0& 0& 0& 0\\
0& 0& 0& 0& 0& 0& 0& 0& 0& 0& 0& -1& 0& 0& 0& 0\\
0& 0& 0& 0& 0& 0& 0& 0& 0& 0& 0& 0& 0& 0& 1& 0\\
0& 0& 0& 0& 0& 0& 0& 0& 0& 0& 0& 0& 0& 0& 0& 1
\end{array}\right).
$$}

Besides, $x_1=x_{\alpha_1}(1)=Q_1\cdot w_1^{-1}$, $w_3=w_{\alpha_3}(1)=w_{1,3}w_1w_{1,3}$, $Q_3=Q_{\alpha_3}=w_{1,3}Q_1w_{1,3}$, $x_{1+2}=x_{\alpha_1+\alpha_2}(1)=w_2x_1w_2^{-1}$, $x_2=x_{\alpha_2}(1)=w_1x_{1+2}w_1^{-1}$.

Let $W_1=(a_{i,j})$, $W_2=(b_{i,j})$. Since $W_1$ commute with $Q_{\alpha_3}$ and $Q_{\alpha_i}$, it follows that $W_1$ on the first basis part is
{\scriptsize $$
\begin{pmatrix}
a_1& a_2& a_3& a_4& a_5-a_1& a_2-a_6& -a_7& a_4-a_8\\
a_9& a_{10}& a_{11}& a_{12}& a_{13}& a_{14}& a_{15}& -a_{16}\\
a_{17}& a_{18}& a_{19}& a_{20}& -a_{21}& a_{18}-a_{22}& a_{23}-a_{19}& a_{20}-a_{24}\\
a_{25}& a_{26}& a_{27}& a_{28}& a_{29}& a_{30}& a_{31}& a_{32}\\
a_1-a_5& a_6& a_7& a_8& a_5& a_2& a_3-a_7& a_4\\
-a_{13}-a_9& -a_{14}& -a_{15}-a_{11}& a_{16}& a_9& a_{10}-a_{14}& a_{11}& a_{12}+a_{16}\\
a_{21}& a_{22}& a_{19}-a_{23}& a_{24}& a_{17}-a_{21}& a_{18}& a_{23}& a_{20}\\
-a_{29}-a_{25}& -a_{30}& -a_{31}-a_{27}& -a_{32}& a_{25}& a_{26}-a_{30}& a_{27}& a_{28}-a_{32}
\end{pmatrix},
$$}
on the second one it is
{\tiny $$
\begin{pmatrix}
-a_{33}+a_{34}& a_{35}+a_{36}& a_{33}& a_{35}& a_{37}& a_{35}+a_{36}-a_{38}& a_{39}& a_{35}+a_{40}-a_{38}\\
-a_{41}-a_{42}& a_{43}+a_{44}& a_{42}& a_{43}& a_{44}& a_{45}& -a_{47}-a_{45}& a_{46}+a_{48}\\
-a_{33}& a_{36}& a_{34}& a_{35}+a_{36}& a_{39}& a_{36}-a_{40}& a_{37}+a_{39}& a_{35}+a_{36}-a_{38}\\
a_{41}& -a_{43}& -a_{41}-a_{42}& a_{44}& a_{47}& -a_{48}-a_{46}& a_{45}& a_{47}\\
-a_{37}& a_{38}& a_{33}-a_{39}& a_{48}& a_{34}+a_{37}-a_{33}&  a_{35}+a_{36}& a_{39}& a_{35}\\
a_{41}+a_{42}-a_{45}& -a_{46}& a_{47}+a_{45}-a_{42}& -a_{48}-a_{46}& -a_{41}-a_{42}& a_{44}-a_{46}& a_{42}& a_{43}-a_{46}-a_{48}\\
a_{39}-a_{33}& a_{40}& a_{33}-a_{37}-a_{39}& a_{38}& -a_{39}& a_{36}& a_{34}+a_{37}+a_{39}-a_{33}& a_{35}+a_{36}\\
-a_{41}-a_{47}& a_{46}+a_{48}& a_{41}+a_{42}-a_{45}& a_{48}& a_{41}& a_{46}+a_{48}-a_{43}& -a_{41}-a_{42}& a_{44}+a_{48}
\end{pmatrix}.
$$}

Let $W_2=(b_{i,j})$, $1\leqslant i,j\leqslant 16$.

Consider the basis change with the matrix~$C$, which on the first part is
{\scriptsize $$
\begin{pmatrix}
1& 0& a_{1,3}& -a_{5,8}& 0& a_{1,10}& a_{1,11}& b_{1,12}\\
a_{1,2}& 1& a_{6,7}& 0& a_{2,9}& 0& a_{6,15}& a_{2,12}\\
-a_{1,3}& a_{5,8}& 1+a_{1,3}& 0& -a_{1,11}& a_{1,10}-b_{1,12}& a_{1,11}& a_{1,10}\\
a_{6,7}-a_{1,2}& 0& a_{1,2}& 1& a_{2,9}+a_{6,15}& a_{2,12}& a_{2,9}& a_{2,12}\\
0& a_{1,10}& a_{1,11}& a_{5,8}+b_{1,12}& 1& -a_{1,3}-a_{1,11}& a_{5,8}& 0\\
a_{1,2}+a_{2,9}& 0& -a_{6,7}+a_{6,15}& a_{2,12}& -a_{1,2}& 1& a_{6,7}& a_{2,12}\\
-a_{1,11}& 0& a_{1,11}& a_{1,10}& a_{1,3}+a_{1,11}& 0& 1& 0\\
0& -a_{2,12}& a_{1,2}+a_{2,9}& -a_{2,12}& 0& -a_{2,12}& -a_{1,2}& 1
\end{pmatrix}
,$$
}
on the second part it is
{\scriptsize $$
\begin{pmatrix}
1& 0& -a_{1,3}& a_{5,8}& 0& -a_{1,10}& -a_{1,11}& -b_{1,12}\\
-a_{1,2}& 1& a_{6,7}& 0& -a_{2,9}& 0& -a_{6,15}& -a_{2,12}\\
a_{1,3}& -a_{5,8}& 1-a_{1,3}& 0& a_{1,11}& -a_{1,10}+b_{1,12}& -a_{1,11}& 0\\
a_{1,2}-a_{6,7}& 0&  a_{1,2}& 1& a_{2,9}+a_{6,15}& a_{2,12}& -a_{2,9}& a_{2,12}\\
0& a_{1,12}& 0& a_{5,8}+b_{1,12}& 1& 0& -a_{1,3}-a_{1,11}& a_{5,8}\\
0& a_{1,2}+a_{2,9}& 0& -a_{6,7}+a_{6,15}& -a_{1,2}& 1& a_{6,7}& a_{2,12}\\
0& 0& 0& 0& 0& 0& 1& 0\\
0& 0& 0& 0& 0& 0& 0& 1
\end{pmatrix}
.$$
}

This matrix commutes with $Q_1$, $Q_3$, $Q_i$ and with all other obtained matrices, so the basis change with it does not move our fixed matrices. From the other side, in linearized form in the matrices $W_1$ and $W_2$ the following elements become zero:
$a_{1,3}$, $a_{1,10}$, $a_{1,11}$, $a_{2,12}$, $b_{1,12}$, $a_{1,2}$, $a_{6,7}$, $a_{2,9}$, $a_{6,15}$, $a_{5,8}$.

Now we introduced how to express the unknown matrices $X_1$, $X_{1+2}$, $X_2$ trough $W_1$, $W_2$:
\begin{align*}
X_1&=Q_1W_1w_1^2+Q_1w_1W_1w_1+Q_1w_1^2W_1;\\
X_{1+2}&=W_2x_1w_2^3+w_2X_1w_2^3+w_2x_1W_2w_2^2+w_2x_1w_2W_2w_2+w_2x_1w_2^2W_2;\\
X_2&=W_1x_{1+2}w_1^3+w_1X_{1+2}w_1^3+w_1x_{1+2}W_1w_1^2+w_1x_{1+2}w_1W_1w_1+w_1x_{1+2}w_1^2W_1.
\end{align*}

Let us write the list of linearized conditions:

$$
\begin{cases}
&W_2Q_i-Q_iW_2=0;\\
&W_1w_2(w_1w_2)^2+w_1W_2(w_1w_2)^2+w_1w_2W_1w_2w_1w_2+w_1w_2w_1W_2w_1w_2+\\
&\quad \quad +(w_1w_2)^2W_1w_2+(w_1w_2)^2w_1W_2=0;\\
&w_{1,3}W_1w_{1,3}w_1+w_3W_1-w_1w_{1,3}W_1w_{1,3}-W_1w_3=0;\\
&w_2w_1w_3W_2+w_2w_1w_{1,3}W_1w_{1,3}w_2+w_2W_1w_3w_2+W_2w_1w_3w_2=0;\\
&x_1X_{1+2}+X_1x_{1+2}-X_{1+2}x_1-x_{1+2}X_1=0;\\
&X_1x_2x_{1+2}+x_1X_2x_{1+2}+x_1x_2X_{1+2}-X_2x_1-x_2X_1=0.
\end{cases}
$$

After that by direct calculus we obtain that the matrices $W_1$ and $W_2$ are zeros. Therefore, $\varphi_2(w_{\alpha_1}(1))=w_{\alpha_1}(1)$ and $\varphi_2(w_{\alpha_2}(1))=w_{\alpha_2}(1)$, what was required.

\subsection{Basis parts of the first form.}

Note that the basis part of the first type for the root system $A_l$ with odd $l$ is just a basis of the root system~$A_3$.
Therefore we consider this system with the assertion that the matrices $Q_1$, $Q_3$, $w_{1,3}$ a mapped to themselves.

So we can suppose that we have three simple roots
$\alpha_1,\alpha_2,\alpha_3$, generating the root system $A_3$,
the weight vector basis is numerated as \begin{multline*}
v_1=v_{\alpha_1}, v_2=v_{-\alpha_1}, v_3=v_{\alpha_2},
v_4=v_{-\alpha_2}, v_5=v_{\alpha_3}, v_6=v_{-\alpha_3},\\
v_7=v_{\alpha_1+\alpha_2}, v_8=v_{-\alpha_1-\alpha_2},
v_9=v_{\alpha_2+\alpha_3}, v_{10}=v_{-\alpha_2-\alpha_3},\\
v_{11}=v_{\alpha_1+\alpha_2+\alpha_3},
v_{12}=v_{-\alpha_1-\alpha_2-\alpha_3}, v_{13}=h_{\alpha_1},
v_{14}=h_{\alpha_2}, v_{15}= h_{\alpha_3}. \end{multline*}

In this basis the matrices representing $w_{\alpha_1}(1)$ and $w_{\alpha_2}(1)$, have the following form

{\scriptsize
$$
w_1=\left(\begin{array}{ccccccccccccccc}
0& -1& 0& 0& 0& 0& 0& 0& 0& 0& 0& 0& 0& 0& 0\\
-1& 0& 0& 0& 0& 0& 0& 0& 0& 0& 0& 0& 0& 0& 0\\
0& 0& 0& 0& 0& 0& 1& 0& 0& 0& 0& 0& 0& 0& 0\\
0& 0& 0& 0& 0& 0& 0& 1& 0& 0& 0& 0& 0& 0& 0\\
0& 0& 0& 0& 1& 0& 0& 0& 0& 0& 0& 0& 0& 0& 0\\
0& 0& 0& 0& 0& 1& 0& 0& 0& 0& 0& 0& 0& 0& 0\\
0& 0& -1& 0& 0& 0& 0& 0& 0& 0& 0& 0& 0& 0& 0\\
0& 0& 0& -1& 0& 0& 0& 0& 0& 0& 0& 0& 0& 0& 0\\
0& 0& 0& 0& 0& 0& 0& 0& 0& 0& 1& 0& 0& 0& 0\\
0& 0& 0& 0& 0& 0& 0& 0& 0& 0& 0& 1& 0& 0& 0\\
0& 0& 0& 0& 0& 0& 0& 0& -1& 0& 0& 0& 0& 0& 0\\
0& 0& 0& 0& 0& 0& 0& 0& 0& -1& 0& 0& 0& 0& 0\\
0& 0& 0& 0& 0& 0& 0& 0& 0& 0& 0& 0& -1& 1& 0\\
0& 0& 0& 0& 0& 0& 0& 0& 0& 0& 0& 0& 0& 1& 0\\
0& 0& 0& 0& 0& 0& 0& 0& 0& 0& 0& 0& 0& 0& 1
\end{array}\right),
$$
}
and
{\scriptsize
$$
w_2=\left(\begin{array}{ccccccccccccccc}
0& 0& 0& 0& 0& 0& 1& 0& 0& 0& 0& 0& 0& 0& 0\\
0& 0& 0& 0& 0& 0& 0& 1& 0& 0& 0& 0& 0& 0& 0\\
0& 0& 0& -1& 0& 0& 0& 0& 0& 0& 0& 0& 0& 0& 0\\
0& 0& -1& 0& 0& 0& 0& 0& 0& 0& 0& 0& 0& 0& 0\\
0& 0& 0& 0& 0& 0& 0& 0& 1& 0& 0& 0& 0& 0& 0\\
0& 0& 0& 0& 0& 0& 0& 0& 0& 1& 0& 0& 0& 0& 0\\
-1& 0& 0& 0& 0& 0& 0& 0& 0& 0& 0& 0& 0& 0& 0\\
0& -1& 0& 0& 0& 0& 0& 0& 0& 0& 0& 0& 0& 0& 0\\
0& 0& 0& 0& -1& 0& 0& 0& 0& 0& 0& 0& 0& 0& 0\\
0& 0& 0& 0& 0& -1& 0& 0& 0& 0& 0& 0& 0& 0& 0\\
0& 0& 0& 0& 0& 0& 0& 0& 0& 0& 1& 0& 0& 0& 0\\
0& 0& 0& 0& 0& 0& 0& 0& 0& 0& 0& 1& 0& 0& 0\\
0& 0& 0& 0& 0& 0& 0& 0& 0& 0& 0& 0& 1& 0& 0\\
0& 0& 0& 0& 0& 0& 0& 0& 0& 0& 0& 0& 1& -1& 1\\
0& 0& 0& 0& 0& 0& 0& 0& 0& 0& 0& 0& 0& 0& 1
\end{array}\right).
$$
}

Besides,  $x_{\alpha_1}(1)$ is
{\scriptsize
$$
x_1=\left(\begin{array}{ccccccccccccccc}
1& -1& 0& 0& 0& 0& 0& 0& 0& 0& 0& 0& -2& 1& 0\\
0& 1& 0& 0& 0& 0& 0& 0& 0& 0& 0& 0& 0& 0& 0\\
0& 0& 1& 0& 0& 0& 0& 0& 0& 0& 0& 0& 0& 0& 0\\
0& 0& 0& 1& 0& 0& 0& 1& 0& 0& 0& 0& 0& 0& 0\\
0& 0& 0& 0& 1& 0& 0& 0& 0& 0& 0& 0& 0& 0& 0\\
0& 0& 0& 0& 0& 1& 0& 0& 0& 0& 0& 0& 0& 0& 0\\
0& 0& -1& 0& 0& 0& 1& 0& 0& 0& 0& 0& 0& 0& 0\\
0& 0& 0& 0& 0& 0& 0& 1& 0& 0& 0& 0& 0& 0& 0\\
0& 0& 0& 0& 0& 0& 0& 0& 1& 0& 0& 0& 0& 0& 0\\
0& 0& 0& 0& 0& 0& 0& 0& 0& 1& 0& 1& 0& 0& 0\\
0& 0& 0& 0& 0& 0& 0& 0& -1& 0& 1& 0& 0& 0& 0\\
0& 0& 0& 0& 0& 0& 0& 0& 0& 0& 0& 1& 0& 0& 0\\
0& 1& 0& 0& 0& 0& 0& 0& 0& 0& 0& 0& 1& 0& 0\\
0& 0& 0& 0& 0& 0& 0& 0& 0& 0& 0& 0& 0& 1& 0\\
0& 0& 0& 0& 0& 0& 0& 0& 0& 0& 0& 0& 0& 0& 1
\end{array}\right);
$$}
the matrix $w_{1,3}$ is
{\scriptsize
$$
\left(\begin{array}{ccccccccccccccc}
0& 0& 0& 0& 1& 0& 0& 0& 0& 0& 0& 0& 0& 0& 0\\
0& 0& 0& 0& 0& 1& 0& 0& 0& 0& 0& 0& 0& 0& 0\\
0& 0& 0& 0& 0& 0& 0& 0& 0& 0& 0& -1& 0& 0& 0\\
0& 0& 0& 0& 0& 0& 0& 0& 0& 0& -1& 0& 0& 0& 0\\
1& 0& 0& 0& 0& 0& 0& 0& 0& 0& 0& 0& 0& 0& 0\\
0& 1& 0& 0& 0& 0& 0& 0& 0& 0& 0& 0& 0& 0& 0\\
0& 0& 0& 0& 0& 0& 0& 1& 0& 0& 0& 0& 0& 0& 0\\
0& 0& 0& 0& 0& 0& 1& 0& 0& 0& 0& 0& 0& 0& 0\\
0& 0& 0& 0& 0& 0& 0& 0& 0& 1& 0& 0& 0& 0& 0\\
0& 0& 0& 0& 0& 0& 0& 0& 1& 0& 0& 0& 0& 0& 0\\
0& 0& 0& -1& 0& 0& 0& 0& 0& 0& 0& 0& 0& 0& 0\\
0& 0& -1& 0& 0& 0& 0& 0& 0& 0& 0& 0& 0& 0& 0\\
0& 0& 0& 0& 0& 0& 0& 0& 0& 0& 0& 0& 0& -1& 1\\
0& 0& 0& 0& 0& 0& 0& 0& 0& 0& 0& 0& 0& -1& 0\\
0& 0& 0& 0& 0& 0& 0& 0& 0& 0& 0& 0& 1& -1& 0
\end{array}\right).
$$}

Consider now the common form of a matrix with entries from the radical, commuting with $Q_1$, $Q_3$ and $w_{1,3}$. By direct calculation we obtain that such a matrix~$C$ is decomposed into two diagonal blocks with respect to the basis parts $\{ v_{\pm 1}, v_{\pm 3}, V_1, V_2,V_3\}$ and $\{ v_{\pm 2}, v_{\pm 4}, v_{\pm 5}, v_{\pm 6}\}$; on the first basis part it is (we do not write coefficients multiplied by two, because they are absent in linearized form)
{\tiny
$$
\begin{pmatrix}
c_{1,1}& c_{1,2}& c_{1,5}& -c_{1,5}& -c_{1,5}& c_{1,14}& c_{1,5}\\
c_{1,2}& c_{1,1} - c_{1,2} + c_{1,5} & -c_{1,5}& c_{1,5}&
-c_{1,5} & -c_{1,2} + c_{1,5}  - c_{1,14}& -c_{1,5}\\
c_{1,5}& -c_{1,5}& c_{1,1}& c_{1,2}& c_{1,5}&  - c_{1,14}& -c_{1,5} \\
-c_{1,5}& c_{1,5}& c_{1,2} & c_{1,1} - c_{1,2} + c_{1,5} & -c_{1,5}&
c_{1,5} - c_{1,2} + c_{1,14}& -c_{1,5} \\
c_{13,1}&-c_{1,2} + c_{13,1} + c_{1,5} & c_{13,5}& -c_{13,5}&
c_{1,1} + c_{1,2} + c_{13,1}& c_{13,14}& c_{13,5}\\
-c_{1,5} & c_{1,5} &  c_{1,5} &  -c1_5 & -c_{1,5}& c_{1,1}+c_{1,2}+c_{1,14}& c_{1,5}\\
-c_{1,5}-c_{13,5}& c_{13,5}+c_{1,5}& c_{1,5}+c_{13,1}& -c_{1,2}+c_{13,1}& -c_{1,5}-c_{13,5}&
c_{15,14}&  c_{15,15}
\end{pmatrix};
$$}
on the second part it is
{\tiny
$$
\begin{pmatrix}
c_{3,3}& c_{3,4}& c_{3,7}& c_{3,8}& c_{3,9}& c_{3,10}& c_{3,11}& c_{3,12}\\
c_{3,4}& c_{3,9} + c_{3,11} + c_{3,3} + c_{3,7}& c_{3,10}-c_{3,4}& c_{3,9}+c_{3,11}& c_{3,8}-c_{3,4}& c_{3,11}+c_{3,7}&
c_{4,10}& c_{3,11}\\
-c_{3,7}& c_{3,4}-c_{3,8}& c_{3,3}+c_{3,7}& c_{3,4}& -c_{3,11}& c_{3,10}-c_{3,12}& c_{3,9}+c_{3,11}& c_{3,10}\\
-c_{3,10}& -c_{3,9}-c_{3,11}& c_{3,4}& c_{3,3}+c_{3,7}& c_{3,10}-c_{3,12}& -c_{3,11}& c_{3,8}-c_{3,4}& c_{3,7}\\
-c_{3,9}& c_{3,4}-c_{3,10}& -c_{3,11}& c_{3,8}-c_{3,12}& c_{3,3}+c_{3,9}& c_{3,4}& c_{3,7}+c_{3,11}& c_{3,8}\\
-c_{3,8}& -c_{3,7}-c_{3,11}& c_{3,8}-c_{3,12}& -c_{3,11}& c_{3,4}& c_{3,3}+c_{3,9}& c_{3,10}-c_{3,4}& c_{3,9}\\
c_{3,11}& c_{3,4}- c_{3,8} + c_{3,12} - c_{3,10}&  -c_{3,9} - c_{3,11}& c_{3,4}-c_{3,10}& -c_{3,7}-c_{3,11}& c_{3,4}-c_{3,8}&
c_{11,11}& c_{3,4}\\
c_{3,12}& c_{3,11}& -c_{3,8}& -c_{3,7}& -c_{3,10}& -c_{3,9}& c_{3,4}& c_{3,3}
\end{pmatrix}.
$$
}

Let
$W_1=\varphi_2(w_{\alpha_1}(1))-w_{\alpha_1}(1)=(x_{i,j})$, $W_2=\varphi_2(w_{\alpha_2}(1))-w_{\alpha_2}(1)=(y_{i,j})$.

The expression of other unknown matrices trough  $W_1$ and $W_2$ is the same as in the previous subsection.

The linearized conditions are:

$$
\begin{cases}
&W_1Q_3-Q_3W_1=0;\\
&W_1w_{\alpha_2}(1)(w_{\alpha_1}(1)w_{\alpha_2}(1))^2+w_{\alpha_1}(1)W_2(w_{\alpha_1}(1)w_{\alpha_2}(1))^2+\\
&\quad \quad+w_{\alpha_1}(1)w_{\alpha_2}(1)W_1w_{\alpha_2}(1)w_{\alpha_1}(1)w_{\alpha_2}(1)+w_{\alpha_1}(1)w_{\alpha_2}(1)w_{\alpha_1}(1)
W_2w_{\alpha_1}(1)w_{\alpha_2}(1)+\\
&\quad \quad +(w_{\alpha_1}(1)w_{\alpha_2}(1))^2W_1w_{\alpha_2}(1)+(w_{\alpha_1}(1)w_{\alpha_2}(1))^2w_{\alpha_1}(1)
W_2=0;\\
&w_{1,3}W_1w_{1,3}w_{\alpha_1}(1)+w_{\alpha_3}(1)W_1-w_{\alpha_1}(1)w_{1,3}W_1w_{1,3}-W_1w_3=0;\\
&w_2w_1w_3W_2+w_2w_1w_{1,3}W_1w_{1,3}w_2+w_2W_1w_3w_2+W_2w_1w_3w_2=0;\\
&x_1X_{1+2}+X_1x_{1+2}-X_{1+2}x_1-x_{1+2}X_1=0;\\
&X_1x_2x_{1+2}+x_1X_2x_{1+2}+x_1x_2X_{1+2}-X_2x_1-x_2X_1=0.
\end{cases}
$$

Let us choose in the matrix $C$ the following coefficients:
$c_{3,3}=-y_{1,7}$; $c_{1,1}=0$; $c_{13,14}=y_{13,13}$;  $c_{3,12}=y_{1,10}$; $c_{13,1}=x_{13,14}$; $c_{13,5}=-x_{15,14}$; $c_{1,5}=-x_{13,5}$; $c_{3,11}=-x_{5,5}$; $c_{3,9}=x_{4,8}$; $c_{3,7}=-x_{3,7}$; $c_{1,2}=-x_{1,14}$; $c_{3,10}=-x_{3,10}$; $c_{3,4}=-x_{3,4}$, after that we conjugate our matrices $w_1+W_1$ and $w_2+W_2$ by the matrix~$E+C$. In the linearized form we obtain
$$
(E+C)(w_i+W_i)(E-C)\sim w_i+W_i+Cw_i-w_iC.
$$

Therefore we see that the matrices $W_1$ and $W_2$ have now the following zero elements (in linearized form):
$y_{5,9}$, $y_{1,7}$, $y_{13,13}$, $y_{1,10}$, $x_{13,14}$,  $x_{15,14}$, $x_{13,5}$, $x_{5,5}$, $x_{4,8}$, $x_{3,7}$, $x_{1,14}$,  $x_{3,10}$, $x_{3,4}$.

After that we directly apply all listed conditions and obtain that matrices $W_1$ and $W_2$ are zeros. Therefore,
$\varphi_2(w_{\alpha_1}(1))=w_{\alpha_1}(1)$ and $\varphi_2(w_{\alpha_2}(1))=w_{\alpha_2}(1)$, what was required.

\section{The images of $x_{\alpha_i}(t)$.}

Now we are interested in the images of the matrices $x_{\alpha_i}(t)$.

Since all basis part are studied similarly, we will consider only the basis part of the second type.

Since all Weil group elements are mapped into themselves under the action of~$\varphi$, it is sufficient to follow the images of
$x_{\alpha_1}(t)$. We fix an arbitrary element $t\in R$ and study the image $\varphi(x_{\alpha_1}(t))=X_t$. Since the matrix~$X_t$ commutes with $w_{\alpha_3}$, $x_{\alpha_1}(1)$, $x_{\alpha_3}(1)$, $x_{\alpha_i}(1)$, $Q_i$, $Q_3$, $x_{\alpha_1+\alpha_2}(1)$ and $x_{-\alpha_2}(1)$, we directly have that $X_t$
is
 {\tiny
$$ \left(\begin{array}{cccccccccccccccc}
x_{1,1}& 0& x_{1,3}& 0& 0& 0& 0& 0& 0& x_{1,10}& 0& 0& 0& 0& 0&0\\
0& x_{2,2}& 0& 0& 0& 0& 0& 0& 0& 0& 0& 0& 0& 0& 0& 0\\
0& 0& x_{1,1}& 0& 0& 0& 0& 0& 0& 0& 0& 0& 0& 0& 0&0\\
0& x_{4,2}& 0& x_{2,2}& 0& 0& 0& 0& 0& 0& x_{4,11}& 0& 0& 0& 0&0\\
0& 0& 0& 0& x_{1,1}& 0& 0&0&  0& 0& 0& 0& 0& 0& 0&0\\
0& 0& 0& 0& 0& x_{2,2}& 0& 0& 0&  0& 0&0& 0&0&0&0\\
0& 0& 0& 0& 0& 0 & x_{1,1}& 0& 0& 0& 0& 0& 0& 0& 0& 0\\
0& 0& 0& 0& 0& 0& 0& x_{2,2}& 0& 0& 0& 0& 0& 0& 0&0\\
0& -x_{1,10}& 0& 0& 0& 0& 0& 0& x_{1,1} &0& x_{1,3} & 0& 0& 0& 0&0\\
0& 0& 0& 0& 0& 0& 0& 0& 0& x_{2,2} & 0& 0& 0& 0& 0& 0\\
0& 0& 0& 0& 0& 0& 0&0&0& 0&  x_{1,1}& 0& 0&0& 0&0\\
0& 0& -x_{4,11}& 0& 0& 0& 0&0& 0& x_{4,2}& 0& x_{2,2}& 0& 0& 0&0\\
0& 0& 0& 0& 0& 0& 0& 0& 0& 0& 0& 0&  x_{1,1}& 0& 0& 0\\
0& 0& 0& 0&0&0&0&0&0&0&0&0&0& x_{2,2}& 0&0\\
0&0&0&0&0& 0& 0&0&0& 0&0& 0& 0&0& x_{1,1}&0\\
0& 0& 0& 0& 0& 0& 0& 0& 0& 0& 0& 0& 0& 0& 0&x_{2,2}
\end{array}\right).
$$
}

Now we introduce $X_t^{1+2}=\varphi_2(x_{\alpha_1+\alpha_2}(t))=w_2 X_t w_2^{-1}$ and take the condition
$X_t^{1+2}x_{\alpha_2}(1)X_t=X_tx_{\alpha_2}(1)$.

Its positions $(1,1)$ and $(2,2)$ imply $x_{1,1}(x_{1,1}-1)=0$ and $x_{2,2}(x_{2,2}-1)=0$, therefore $x_{1,1}=x_{2,2}=1$. From the position $(1,10)$ it follows $x_{1,10}=0$, from the position $(4,13)$ it follows $x_{4,11}=0$. Besides, commuting with other Weil group elements gives $x_{4_2}=-x_{1,3}$, so $X_t=x_{\alpha_1}(s)$ for some $s\in R^*$.

If we introduce the image of $h_{\alpha_1}(t)$, then by the same arguments we obtain that it is $h_{\alpha_1}(s)$.

 Studying other basis parts we can easily check that for every $t\in R^*$ the corresponding $s\in R^*$ is unique for the whole basis.

\section{Proof of the main theorem.}

Clear that $\varphi_2(h_{\alpha_k}(t))=h_{\alpha_k}(s)$,
$k=1,\dots,n$. Denote the mapping $t\mapsto s$ by $\rho: R^*
\to R^*$. Note that for $t\in R^*$ we have
$\varphi_2(x_1(t))=\varphi_2(h_{\alpha_2}(t^{-1}) x_1(1)
h_{\alpha_2}(t))=h_{\alpha_2}(s^{-1}) x_1(1) h_{\alpha_2}(s)=
x_1(s)$. If $t\notin R^*$, then $t\in J$, i.\,e., $t=1+t_1$, where
$t_1\in R^*$. Then
$\varphi_2(x_1(t))=\varphi_2(x_1(1)x_1(t_1))=x_1(1)x_1(\rho(t_1))=
x_1(1+\rho(t_1))$. Therefore if we extend the mapping $\rho$ on the whole ring~$R$ (by the formula $\rho(t):=1+\rho(t-1)$, $t\in
R$), we obtain $\varphi_2(x_1(t))=x_1(\rho(t))$ for all $t\in R$.
Clear that $\rho$ is injective, additive and also multiplicative on all invertible elements. Since every element of~$R$ is a sum of two invertible elements, we see that $\rho$ is an isomorphism from the ring~$R$ onto some its subring~$R'$. Note that in this situation $C E(\Phi,R) C^{-1}=E(\Phi,R')$ for some matrix $C\in
\GL(V)$. Show that $R'=R$.

Let us denote the matrix units by~$E_{ij}$.

\begin{lemma}\label{Tema}
If for some $C\in \GL(V)$ we have $C E(\Phi,R) C^{-1}=
E(\Phi,R')$, where $R'$ is a subring in~$R$, then $R'=R$.
\end{lemma}
\begin{proof}
Show that with the group $E_{\ad}(\Phi,R)$ of the case under consideration by addition and multiplication of matrices we can we can obtain every matrix $a
E_{1,1}$, $a\in R$.

The matrix $((x_{\alpha_1}(1)-1)(x_{\alpha_2}-1))^2$ has the single nonzero element $\cdot E_{7,8}$. Similarly,
$E_{8,7}=((x_{-\alpha_1}-1)(x_{-\alpha_2}-1))^2$. Therefore,
we directly obtain the matrix units $E_{7,7}, E_{7,8}, E_{8,7}$,
$E_{8,8}$. According to the transitive action of the Weil group on the roots,
by conjugating the given matrix units by different Weil group elements, we obtain $E_{i,j}$, $i,j\in \{ 1,2\}$, $\{ 3,4\}$, $\{
5,6\}$, $\{ 9,10\}$, or $\{ 11,12\}$. So we have $E_{1,1}$. Later with the help of  $E_{1,1} h_{\alpha_2}(t) E_{1,1}$ we get $t
E_{1,1}$ for all invertible~$t\in R$. Since every element of~$R$ is a sum of two invertible elements, we obtain all $aE_{1,1}$, $a\in R$.

Suppose now that $R'$ is a proper subring of~$R$.

Since $C E(\Phi,R) C^{-1}= E(\Phi,R')$, then the subring of
$M_n(R)$, generated by all elements of $E(\Phi,R)$, is mapped to the subring of $M_n(R')$, generated by all elements of $E(\Phi,R')$. Therefore all coefficients of
$G_a=C (aE_{1,1}) C^{-1}$ lie in the subring~$R'$.

Let $C=(c_{k,l})$, $C^{-1}=(c_{k,l}')$. Note that the  $(i,j)$-th
coefficient of $G_a$ is $ac_{i,1}c_{1,j}$.  Every invertible matrix over a local ring in every its row and every its line has at least one invertible coefficient, so there exist such indices  $i_0$ and
$j_0$, for which $c_{i_0,1}c_{1,j_0}$ is invertible. On the corresponding place in the matrix $G_a$ it is an arbitrary element of~$R$.

The obtained contradiction shows that $R'=R$.
\end{proof}

Therefore we proved that   $\rho$ is an automorphism of~$R$. Consequently the composition of the initial automorphism~$\varphi$ and some basis change with the help of a matrix $C\in \GL_n(R)$, (mapping $E(\Phi,R)$ into itself)
is a ring automorphism~$\rho$. It proves Theorem~2. $\square$

\bigskip

The main theorem (Theorem 1) now obviously follows from Theorem~2, proved above, and also Theorems~1 and~3, proved in the paper~\cite{normalizers}. Since in~\cite{normalizers} in the proof of Theorems~1 and~3 we did not used that two was invertible, and only required Theorem~2 (similar to our Theorem~2) of the given paper, then the proof is now automatically extended to the case of rings without~$1/2$.

\end{document}